\documentclass[12pt]{article}
\usepackage{amsmath,amsfonts,amssymb,mathtools}
\usepackage{amsthm}
\usepackage{amsmath}
\usepackage{authblk}
\usepackage[hmargin=2.5cm,vmargin=3cm,centering]{geometry}
\usepackage{hyperref}
\usepackage{color}
\usepackage{cite}
\allowdisplaybreaks

\hypersetup{
    colorlinks,
    citecolor=black,
    filecolor=magenta,
    linkcolor=black,
    urlcolor=black,
    }
\usepackage{color}
 \definecolor{red}{rgb}{1,0,0}

\newcommand{\xddots}{%
  \raise 4pt \hbox {.}
  \mkern 6mu
  \raise 1pt \hbox {.}
  \mkern 6mu
  \raise -2pt \hbox {.}
}

\theoremstyle{break}
\newtheorem{thm}{Theorem}[section] 
 
\newtheorem{RQS}{Remarks}[thm]
\newtheorem{NTS}{Notations}[thm]
 
\newtheorem{defi}[thm]{Definition}

\newtheorem{lem}[thm]{Lemma}
\numberwithin{equation}{section}

\title{Unique weak solutions of the d-dimensional micropolar equation with fractional dissipation}

\title{\large\bf{Unique weak solutions of the d-dimensional micropolar equation with fractional
dissipation}}
\usepackage{authblk}

\author[1]{Oussama~Ben~Said \thanks{obensai@okstate.edu}}
\author[2]{Jiahong Wu \thanks{jiahong.wu@okstate.edu }}
\affil[1,2]{Department of Mathematics, 
	Oklahoma State University,
	Stillwater, OK 74078, USA}
\usepackage{titlesec} 

\titleformat{\subsubsection}[runin]
  {\normalfont\large\bfseries}{\thesubsubsection}{1em}{}
\begin{document}

\maketitle
\begin{abstract}
This article examines the existence and uniqueness of weak solutions to the d-dimensional micropolar equations ($d=2$ or $d=3$) with general fractional dissipation $(-\Delta)^{\alpha}u$ and $(-\Delta)^{\beta}w$. The micropolar equations with standard Laplacian dissipation model
fluids with microstructure. The generalization to include fractional dissipation allows simultaneous study of a family of equations and is relevant in some physical circumstances. We establish that, when $\alpha\ge \frac12$ and $\beta\ge \frac12$, any initial data $(u_0, w_0)$ 
in the critical Besov space $u_0\in B^{1+\frac{d}{2}-2\alpha}_{2,1}(\mathbb R^d)$ and $w_0\in B^{1+\frac{d}{2}-2\beta}_{2,1}(\mathbb R^d)$ yields a unique weak solution. For $\alpha\ge 1$ and $\beta=0$, any initial data $u_0\in B^{1+\frac{d}{2}-2\alpha}_{2,1}(\mathbb R^d)$ and $w_0\in B^{\frac{d}{2}}_{2,1}(\mathbb R^d)$ also leads to a unique weak solution as well. The regularity indices in these Besov spaces appear to be optimal and can not be lowered in order to achieve the uniqueness.  Especially,  the 2D micropolar equations with the standard Laplacian dissipation, namely $\alpha=\beta=1$ have a unique weak solution for $(u_0, w_0)\in B^0_{2,1}$. The proof involves the construction of 
successive approximation sequences and extensive {\it a priori} estimates in Besov space settings. 
\end{abstract}

\noindent\textbf{Key words:}\;Micropolar equations, Littlewood-Paley, Local solution,\;Uniqueness.\\
\noindent\textbf{MSC-2010:} 35A01, 35A02, 35Q35, 76D03

\section{Introduction}
\label{int}
\setcounter{equation}{0}

The micropolar equations were first proposed in 1965 by C.A Eringen to modal micropolar fluids which are fluids with microstructure (see \cite{Erdo,Eringen1,Er2,Cowin}). 
These equations can model large number of complex fluids such as  animal blood, suspensions, and liquid crystals.  In this paper, we focus on the following d-dimensional ($d = 2$ or $d = 3$) incompressible micropolar equations with fractional dissipation
\begin{align}
 \left\{
    \begin{array}{ll}
       \partial_tu+(\nu+k)(-\Delta)^{\alpha}u+u\cdot \nabla u+\nabla\Pi-2k\nabla\times w=0,\\\\
         \partial_tw+4kw +2k\nabla\times u+u\cdot \nabla w+\gamma(-\Delta)^{\beta}w =0,\\\\ \nabla\cdot u=0,\\\\
         u(x,0)=u_0(x)\,,\;\;w(x,0)=w_0(x), 
 \end{array} \label{2.1}
\right.
\end{align}
where $u = u(x,t)\in\mathbb{R}^d$ denotes the fluid velocity, $w = w(x,t)\in\mathbb{R}^d$ the field of microrotation representing the angular velocity of the rotation of the fluid particles, $\Pi=\Pi(x,t)$ the scalar pressure, and the parameter $\nu$ denotes the Newtonian kinematic viscosity, $k$ the microrotation viscosity and $\gamma$ the angular viscosity. Here the fractional Laplacian operator $(-\Delta)^{\alpha}$ (which is also referred as the Riesz  potential operator) is defined via the Fourrier transform 
\begin{align*}
\widehat{(-\Delta)^{\alpha}}f(\xi)=|\xi|^{2\alpha}\widehat{f}(\xi)\;,
\end{align*}
where 
\begin{align*}
\widehat{f}(\xi)=\frac{1}{(2\pi)^{\frac{d}{2}}}\int_{\mathbb{R}^d}e^{-ix\cdot\xi}f(x)\,dx\;.
\end{align*}

\vskip .1in
Besides their many physical applications, the micropolar equations are also of great interest 
in mathematics. Fundamental issues such as the well-posedness problem on (\ref{2.1}) 
have recently attracted considerable interest and an array of important results have been 
established (see, e.g., \cite{Boldrinidr,ChenChen,Chenmm,EGaldir,Lukaszewicz1, Lukaszewicz2,Lukaszewicz3,Rojaszammmc,Nowakowski,Yamaguchi}). More recent focuses have been on the micropolar equations with partial or fractional dissipation (see, e.g., \cite{DLW,Dwxye,DZ,JiuLiu,LiuWang,Xue}). Investigations on nonlocal diffusion have now become a trend \cite{BValdinoci}. The study of fractionally dissipated micropolar equations allow us to simultaneously treat a family of equations including those with the standard Laplacian dissipation. The investigations
on the micropolar 
equations with fractional dissipation help reveal how the global well-posedness problem is related 
to the fractional regularization. 

\vskip .1in 
The main goal of this study is to obtain the existence and uniqueness of solutions 
to (\ref{2.1}) in a weakest possible functional setting for the largest possible ranges of $\alpha$ and $\beta$. Our main results can be stated as follows.

\begin{thm}\label{thm}

Consider (\ref{2.1}) with $\alpha \ge \frac{1}{2}$ and $\beta \ge \frac12$. Assume the initial data $u_0$ and $w_0$ satisfy
$$ \nabla\cdot u_0=0, \quad u_0\in B^{\frac{d}{2}+1-2\alpha}_{2,1}(\mathbb R^d), \quad w_0\in B^{\frac{d}{2}+1-2\beta}_{2,1}(\mathbb R^d).
$$
Then there exist $T >0$ and a unique weak solution $(u,w)$ of (\ref{2.1}) on $\left[ 0,T\right]$ satisfying
\begin{align}
&u\in L^{\infty}(0,T,B_{2,1}^{1+\frac{d}{2}-2\alpha}(\mathbb R^d))\cap L^{1}(0,T,B_{2,1}^{1+\frac{d}{2}}(\mathbb R^d))\;,
\label{2.2}\\&
w\in L^{\infty}(0,T,B_{2,1}^{1+\frac{d}{2}-2\beta}(\mathbb R^d))\cap L^{1}(0,T,B_{2,1}^{1+\frac{d}{2}}(\mathbb R^d)).
\label{2.3}\end{align}
\end{thm}

Here $B^r_{p,q}$ denotes the inhomogeneous Besov space. A review of the Besov spaces and 
related facts is provided in the following section. As a special consequence of 
Theorem \ref{thm}, the two-dimensional (2D) micropolar equations with $\alpha=\beta=1$, 
namely the standard Laplacian dissipation always possess a unique local solution $(u, w)$ in the 
critical Besov space $L^\infty(0, T; B^{0}_{2,1} (\mathbb R^2))$. For the 3D micropolar equations
with the standard Laplacian dissipation, the uniqueness is also attained in the critical
Besov space $L^\infty(0, T; B^{\frac12}_{2,1} (\mathbb R^3))$. Here the critical Besov spaces 
are the Besov space settings for which the solution of the differential equations and its scaling 
invariant counterparts share the same norm. In the general fractional dissipation cases, the regularity indices $1+\frac{d}{2}- 2\alpha$ and $1+\frac{d}{2}- 2\beta$ in the Besov spaces appear to be optimal and one may not be able to achieve the uniqueness when they are lowered.

\vskip .1in 
\begin{thm}\label{thm2}
Consider (\ref{2.1}) with  $\alpha \ge 1$ and $\beta = 0$. Assume the initial data $(u_0, w_0)$ satisfies 
$$
\nabla\cdot u_0=0, \quad u_0\in B^{1+\frac{d}{2}-2\alpha}_{2,1}(\mathbb R^d), \quad w_0\in B^{\frac{d}{2}}_{2,1}(\mathbb R^d).
$$
Then there exist $T >0$ and a unique weak solution $(u,w)$ of (\ref{2.1})  on $\left[ 0,T\right]$ satisfying
\begin{align}
&u\in L^{\infty}(0,T,B_{2,1}^{1+\frac{d}{2}-2\alpha}(\mathbb R^d))\cap L^{1}(0,T,B_{2,1}^{1+\frac{d}{2}}(\mathbb R^d))\;,
\label{2.2}\\&
w\in L^{\infty}(0,T,B_{2,1}^{\frac{d}{2}}(\mathbb R^d)).
\label{2.3}\end{align}
\end{thm}

Theorem \ref{thm2} deals with the case when the equation of $w$ involves no diffusion. The Besov space for $u$ remains the same, but the setting for $w$ needs to be in a more regular Besov space due to the lack of diffusion in the equation of $w$. For an inviscid equation, the regularity index $\frac{d}{2}$ in the Besov space $B^{\frac{d}{2}}_{2,1}(\mathbb R^d)$  can not be lowered in order to obtain the uniqueness of solutions. 

\vskip .1in 
The proof for each of the theorems is naturally split into two parts: the existence and uniqueness parts. The existence part starts with the construction a successive approximation 
sequence which iteratively solves systems close to (\ref{2.1}). This successive approximation sequence is then shown to be uniformly bounded in suitable Besov spaces via the method of mathematical induction. These bounds allow us to extract a subsequence, which converges weakly to a limit. Using the Aubin-Lions Lemma, the weak limit is then shown to be the weak solution of (\ref{2.1}). The main efforts are devoted to proving the uniform boundedness. This process 
involves various analysis tools and techniques. The uniqueness is established by analyzing the differences 
in the $L^2$ space.

\vskip .1in 
The rest of this paper is divided into three sections. The second section serves as 
a preparation. It reviews the Besov space and related tools to be used in the subsequent sections. The third section proves Theorem \ref{thm}. It is further divided into two subsections with one devoted to the existence and the other to the uniqueness.  The last section provides the 
proof of Theorem \ref{thm2}. It is again split into two subsections, one for the proof of existence and one for the uniqueness.

\vskip .3in 
\section{Preparations: Besov spaces}

This section serves as a preparation. Materials presented here will be used in 
the proofs of Theorems \ref{thm} and \ref{thm2}. The definition of the Besov space and 
related simple facts can be found in \cite{BCD}. Lemma \ref{lem2.3} is taken from \cite[Lemma A.5] {JiuSuo}. In what follows, $\mathcal S(\mathbb{R}^d)$ denotes the Schwartz class and $\mathcal{S}'(\mathbb{R}^d)$ the tempered distribution.

\begin{defi}[Inhomogenous Besov space $B^s_{p,q}$]

	 $f\in \mathcal{S}'(\mathbb{R}^d)$ belongs to $B^s_{p,q}$ with $s\in \mathbb{R}$ and $1\le p\le q\le \infty$ if 
	\begin{align*}
	\| f\|_{B^s_{p,q}}\equiv\| 2^{sj}\|\Delta_j f\|_{L^p}\|_{l^q}=\left\{
	\begin{array}{ll}
	\Big(\sum\limits_{j=-1}^{+\infty}( 2^{sj}\|\Delta_j f\|_{L^p})^q \Big)^{\frac{1}{q}}\quad & \mbox{\text{if}}\quad q<\infty\;,\\
	\sup\limits_{j\ge -1}\,2^{sj}\|\Delta_j f\|_{L^p}\quad & \mbox{if }\quad q=\infty
	\end{array}
	\right.
	\end{align*}
	is finite.
\end{defi}

\begin{lem}\label{lemma1}
	Let $B(0,r)$ and $C(0,r_1,r_2)$ denote the standard ball and the annulus,  respectively,
	\begin{align*}
	B(0,r)=\left\{\xi\in \mathbb{R}^d,\;\;|\xi|\le r\right\}\;,\quad C(0,r_1,r_2)=\left\{\xi\in \mathbb{R}^d,\;\;r_1\le|\xi|\le r_2\right\}.
	\end{align*}
	There are two compactly supported smooth radial functions $\phi$ and $\psi$ satisfying 
	\begin{align}
	&\mathrm{supp}\,\phi\subset B(0,\frac{4}{3})\;, \quad \mathrm{supp}\,\psi\subset C(0,\frac{3}{4},\frac{8}{3})\;,
	\nonumber\\&
	\phi(\xi)+\sum_{j\ge0}\psi(2^{-j}\xi)=1\quad\text{for all}\;\;\xi\in\mathbb{R}^d\;. \label{1t}
	\end{align}
\end{lem}
The proof of Lemma \ref{lemma1} can be found in \cite[p.59]{BCD}.

\begin{NTS}
	We use $\widetilde{h}$ and $h$ to denote  the inverse Fourier transforms of $\phi$ and $\psi$ respectively 
	\begin{align*}
	\widetilde{h}=\mathcal{F}^{-1}\phi\;,\quad h=\mathcal{F}^{-1}\psi\;.
	\end{align*}
	We write $\psi_j(\xi)=\psi(2^{-j}\xi).$ By a simple property of the Fourier transform, 
	\begin{align*}
	h_j(x):=\mathcal{F}^{-1}\psi_j(x)=2^{dj}h(2^{j}x)\;.
	\end{align*}
\end{NTS}

\begin{defi}
	
	The inhomogeneous dyadic block operator $\Delta_j$ are defined as 
	
	\begin{align*}
	&\Delta_jf=0\,\,\quad\quad\quad\quad\quad\quad\quad\quad\quad\quad\quad\quad\quad\quad\quad\quad\text{for} \; j \le -2 \;,
	\\
	&\Delta_{-1}f=\widetilde{h}\ast f=\int_{\mathbb{R}^d}f(x-y)\widetilde{h}(y)\,dy\,,
	\\
	&\Delta_{j}f=h_j\ast f=2^{dj
	}\int_{\mathbb{R}^d}f(x-y)h(2^{j
	}y)\,dy\quad\quad\;\;\text{for} \; j \ge 0 \;.
	\end{align*}
	The corresponding inhomogeneous low frequency cut-off operator $S_j$ is defined by 
	\begin{align*}
	S_jf=\sum_{k\le j-1} \Delta_{k}f\;.
	\end{align*}
\end{defi}

\begin{RQS}For any function $f$ in the usual Schwarz class $\mathcal{S},$ (\ref{1t}) implies 
	\begin{align*}
	\widehat{f}(\xi)=\phi(\xi)\widehat{f}(\xi)+\sum_{j\ge0}\psi(2^{-j}\xi)\widehat{f}(\xi)\;,
	\end{align*}
	or in terms of the inhomogenous dyadic block operators 
	\begin{align*}
	f=\sum_{j\ge-1}\Delta_j f \quad \text{or}\quad Id=\sum_{j\ge-1}\Delta_j\;,
	\end{align*}
	where $Id$ denotes the identity operator.  For generality, for any $F$ in the space of tempered distributions $\mathcal{S}',$
	\begin{align}
	F=\sum_{j\ge-1}\Delta_j F \quad \text{or}\quad Id=\sum_{j\ge-1}\Delta_j\quad \text{in}\quad\mathcal{S}'\;.
	\label{600}\end{align}
	$(\ref{600})$ is referred to as the Littlewood-Paley decomposition for tempered distributions.
\end{RQS}

\begin{defi}
	In terms of inhomogeneous dyadic block operators, we can write the standard product in terms of the paraproducts, namely
	\begin{align*}
	FG=\sum_{|j-k|< 2}S_{k-1}F\Delta_kG+\sum_{|j-k|< 2}\Delta_kFS_{k-1}G+\sum_{k\ge j-1}\Delta_kF\widetilde{\Delta}_kG\;,
	\end{align*}
	where $\widetilde{\Delta}_k=\Delta_{k-1}+\Delta_{k}+\Delta_{k+1}.$  This is the so-called Bony decomposition. 
\end{defi}

\begin{lem}
	Let $\alpha\ge 0\,.$ Let $1\le p\le q \le\infty\,.$
	\\
	(1) If $f$ satisfies 
	\begin{align*}
	\mathrm{supp}\, \widetilde{f}\subset \left\{\xi \in \mathbb{R}^d,\quad |\xi|\le K2^j\right\}\;,
	\end{align*}
	for some integer $j$ and a constant $K >0$ then  
	\begin{align*}
	\|(-\Delta)^\alpha f\|_{L^q(\mathbb{R}^d)}\le c_12^{2\alpha j+jd(\frac{1}{p}-\frac{1}{q})}\| f\|_{L^p(\mathbb{R}^d)}\;.
	\end{align*} 
	(2) If $f$ satisfies 
	\begin{align*}
	\mathrm{supp}\, \widetilde{f}\subset \left\{\xi \in \mathbb{R}^d,\quad  K_12^j\le|\xi|\le K_22^j\right\}\;,
	\end{align*}
	for some integer $j$ and constants $0<K_1\le K_2$ then  
	\begin{align*}
	c_12^{2\alpha j}\| f\|_{L^q(\mathbb{R}^d)}\le \|(-\Delta)^\alpha f\|_{L^q(\mathbb{R}^d)}\le c_22^{2\alpha j+jd(\frac{1}{p}-\frac{1}{q})}\| f\|_{L^p(\mathbb{R}^d)}\;,
	\end{align*}
	where $c_1,c_2$ are constants depending only on $\alpha, p, q.$
\end{lem}

\vskip .1in 
Below we state bounds for the triple products involving Fourier localized functions. These bounds will be used in the proofs of Theorems \ref{thm} and \ref{thm2}. We refer the reader to Lemma A.5 in \cite{JiuSuo} for a detailed proof of the following lemma.
\begin{lem}\label{lem2.3}
	Let $j\ge 0$ be an integer. Let $\Delta_j$ be the inhomogeneous Littlewood-Paley-localization operator. For any  vectors field $F, G,H$ with $\nabla\cdot F=0$ we have 
	\begin{align*}
	|\int_{\mathbb{R}^d}&\Delta_j(F\cdot \nabla G)\cdot\Delta_jH\, dx|\le c\| \Delta_j H\|_{L^2}\Big(2^j\sum_{m\le j-1}2^{\frac{d}{2}m}\| \Delta_m F\|_{L^2}\sum_{|j-k|\le 2 }\| \Delta_k G\|_{L^2}\\&+\sum_{|j-k|\le 2 }\| \Delta_k F\|_{L^2}\sum_{m\le j-1}2^{(1+\frac{d}{2})m}\| \Delta_m G\|_{L^2}+\sum_{k\le j-1}2^j2^{\frac{d}{2}k}\| \Delta_k F\|_{L^2}\| \widetilde{\Delta}_k G\|_{L^2}\Big)
	\end{align*}
	
	and 
	
	\begin{align*}
	|\int_{\mathbb{R}^d}&\Delta_j(F\cdot \nabla G)\cdot\Delta_jG\, dx|\le c\| \Delta_j G\|_{L^2}\Big(\sum_{m\le j-1}2^{(1+\frac{d}{2})m}\| \Delta_m F\|_{L^2}\sum_{|j-k|\le 2 }\| \Delta_k G\|_{L^2}\\&+\sum_{|j-k|\le 2 }\| \Delta_k F\|_{L^2}\sum_{m\le j-1}2^{(1+\frac{d}{2})m}\| \Delta_m G\|_{L^2}+\sum_{k\le j-1}2^j2^{\frac{d}{2}k}\| \Delta_k F\|_{L^2}\| \widetilde{\Delta}_k G\|_{L^2}\Big)~.
	\end{align*}
\end{lem}

\vskip .3in 
\section{Proof of Theorem \ref{thm}}

\subsection{Existence of a weak solution}

This subsection proves the existence part of Theorem \ref{thm}. The approach is to construct a successive approximation sequence and show that the limit of a subsequence actually solves (\ref{2.1}) in the weak sense.
\begin{proof}[Proof for the existence part of Theorem \ref{thm}] We consider a successive approximation $\left\{(u^{(n)},w^{(n)})\right\}$ satisfying 
\begin{align}
 \left\{
  \begin{array}{ll}
 u^{(1)}=S_2u_0\;,\quad \quad w^{(1)}=S_2w_0\;,\\\\
  \partial_tu^{(n+1)}+(\nu+k)(-\Delta)^{\alpha}u^{(n+1)}=\mathbb{P}(-u^{(n)}\cdot\nabla u^{(n+1)})+2k\nabla \times w^{(n)}\;,\\\\
 \partial_tw^{(n+1)}+\gamma (-\Delta)^{\beta}w^{(n+1)}=-4k w^{(n+1)}-2k\nabla\times u^{(n)}-u^{(n)}\cdot\nabla w^{(n+1)}\;,\\\\
 u^{(n+1)}(x,0)=S_{n+1}u_0\;,\quad\quad w^{(n+1)}(x,0)=S_{n+1}w_0~,
 \end{array}\label{2.6}
\right.
\end{align}
where $\mathbb{P}=I-\nabla(-\Delta)^{-1}\mathrm{div}$ is the standard Leray Projection. For \begin{align*}M=2(\|u_0\|_{B_{2,1}^{1+\frac{d}{2}-2\alpha}}+\|w_0\|_{B_{2,1}^{1+\frac{d}{2}-2\beta}})\;,\end{align*}$T>0$ sufficiently small and $0<\delta<1$  (to be specified later), we set 
\begin{align}
Y\equiv&\left\{(u,w)\;|\quad  \|u\|_{L^{\infty}(0,T,B_{2,1}^{1+\frac{d}{2}-2\alpha})}\le M\;,\;\;\|w\|_{L^{\infty}(0,T,B_{2,1}^{1+\frac{d}{2}-2\beta})}\le M\;,\right.\nonumber \\
&\hphantom{M(u,v)\;|}\left.\quad\|u\|_{L^{1}(0,T,B_{2,1}^{1+\frac{d}{2}}) }\le \delta\;,\;\;\|w\|_{L^{1}(0,T,B_{2,1}^{1+\frac{d}{2}}) }\le \delta\right\}.\label{2.7}
\end{align}
We show that  $\left\{(u^{(n)},w^{(n)})\right\}$ has a subsequence that converges to the weak solution of (\ref{2.1}). This process consists of three main steps. The first step is to show that $\left\{(u^{(n)},w^{(n)})\right\}$ is uniformly bounded in $Y$. The second step is to extract a strongly convergent subsequence via the Aubin-Lions Lemma. While the last step is to show that the limit is indeed a weak solution of (\ref{2.1}).

\vskip .1in
To show the uniform bound for $\left\{(u^{(n)},w^{(n)})\right\}$ in Y, we prove by induction. Clearly,
\begin{align*}
&\|u^{(1)}\|_{L^{\infty}(0,T,B_{2,1}^{1+\frac{d}{2}-2\alpha} )}=\|S_2u_0\|_{L^{\infty}(0,T,B_{2,1}^{1+\frac{d}{2}-2\alpha} )}\le M\;,\\&\|w^{(1)}\|_{L^{\infty}(0,T,B_{2,1}^{1+\frac{d}{2}-2\beta} )}=\|S_2w_0\|_{L^{\infty}(0,T,B_{2,1}^{1+\frac{d}{2}-2\beta} )}\le M\;.
\end{align*}
If $T>0$ is sufficiently small, then 
\begin{align*}
\|u^{(1)}\|_{L^1(0,T,B_{2,1}^{1+\frac{d}{2}} )}\le T \|S_2u_0\|_{B_{2,1}^{1+\frac{d}{2}}}\le T\,c \,\|u_0\|_{B_{2,1}^{1+\frac{d}{2}-2\alpha}}\le \delta\;,
\end{align*}
\begin{align*}
\|w^{(1)}\|_{L^1(0,T,B_{2,1}^{1+\frac{d}{2}} )}\le T \|S_2w_0\|_{B_{2,1}^{1+\frac{d}{2}}}\le T\,c\,\|w_0\|_{B_{2,1}^{1+\frac{d}{2}-2\beta}}\le \delta~.
\end{align*}
Assuming that $(u^{(n)},w^{(n)})$ obeys the bounds defined in $Y,$ namely 

\begin{align*}
&\|u^{(n)}\|_{L^{\infty}(0,T, B_{2,1}^{1+\frac{d}{2}-2\alpha} )}\le M\;,\;\;\|w^{(1)}\|_{L^{\infty}(0,T,B_{2,1}^{1+\frac{d}{2}-2\beta} )}\le M\;,\\&\|u^{(n)}\|_{L^{1}(0,T,B_{2,1}^{1+\frac{d}{2}}) }\le \delta\;,\;\;\|w^{(n)}\|_{L^{1}(0,T,B_{2,1}^{1+\frac{d}{2}}) }\le \delta\;,
\end{align*}
we prove that $\left\{(u^{(n+1)},w^{(n+1)})\right\}$ obeys the same bound for suitably selected $T>0,\; M>0$ and $\delta>0$. For the sake of clarity, the proof of the four bounds in achieved in the following four steps.
\subsubsection{\textbf{The estimate of $u^{(n+1)}$ in $B_{2,1}^{1+\frac{d}{2}-2\alpha}(\mathbb{R}^d)\,.$}} Let $j\ge 0$ be an integer. Applying $\Delta_j$ to the second equation in (\ref{2.6}) and then dotting with $\Delta_ju^{(n+1)},$ we obtain 
\begin{align}
\frac{1}{2}\frac{d}{dt}\|\Delta_ju^{(n+1)}\|^2_{L^2}+(\nu+k)\|\Lambda^{\alpha}\Delta u^{(n+1)}\|^2_{L^2}=A_1+A_2\;,\label{2.8}
\end{align}
where \begin{align*}
&A_1=\int_{\mathbb{R}^d}2k\Delta_j(\nabla \times w^{(n)})\cdot\Delta_j u^{(n+1)}\,dx\;,\\&A_2=-\int\Delta_j(u^{(n)}\cdot\nabla u^{(n+1)})\Delta_ju^{(n+1)}\,dx\;.
\end{align*}
We remark that the projection operator $\mathbb{P}$ has been eliminated due to the divergence-free condition $\nabla\cdot u^{(n+1)}=0.$ The dissipative part admits a lower bound 
\begin{align*}
(\nu + k) \|\Lambda^{\alpha}\Delta_ju^{(n+1)}\|^2_{L^2}\ge c_0\,2^{2\alpha j}\|\Delta_j u^{(n+1)}\|^2_{L^2}\;,
\end{align*} 
where $c_0>0$ is a constant. By H\"older's inequality and Bernstein's inequality
\begin{align*}
|A_1|&=|\int_{\mathbb{R}^d}2k\,\Delta_j(\nabla\times w^{(n)})\cdot \Delta_j u^{(n+1)}\, dx|\\&\le 2k\,\|\Delta_j(\nabla\times w^{(n)})\|_{L^2}\|\Delta_ju^{(n+1)}\|_{L^2}\\&\le c\,2^j\,\|\Delta_jw^{(n)}\|_{L^2}\|\Delta_j u^{(n+1)}\|_{L^2}~.
\end{align*}
According to Lemma \ref{lem2.3}, 
\begin{align*}
|A_2|&=|-\int_{\mathbb{R}^d}\Delta_j(u^{(n)}\cdot \nabla u^{(n+1)})\cdot \Delta_j u^{(n+1)}\,dx|\\&\le c\,\|\Delta_j u^{(n+1)}\|^2_{L^2}\sum_{m\le j-1}2^{(1+\frac{d}{2})m}\|\Delta_m u^{(n)}\|_{L^2}\\&\quad+c\,\|\Delta_ju^{(n+1)}\|_{L^2}\|\Delta_ju^{(n)}\|_{L^2}\sum_{m\le j}2^{(1+\frac{d}{2})m}\|\Delta_m u^{(n+1)}\|_{L^2}\\&\quad+c\sum_{k\ge j-1} 2^j2^{\frac{d}{2}k}\|\Delta_ku^{(n)}\|_{L^2}\|\widetilde{\Delta}_ku^{(n+1)}\|_{L^2}\|\Delta_j u^{(n+1)}\|_{L^2}\;.
\end{align*}
Inserting the estimates above in (\ref{2.8}) and eliminating $\|\Delta_j u^{(n+1)}\|_{L^2}$ from the both sides\;, we obtain 
\begin{align}
\frac{d}{dt}\|\Delta_j u^{(n+1)}\|_{L^2}+c_02^{2\alpha j}\|\Delta_j u^{(n+1)}\|_{L^2}\le J_1+J_2+J_3+J_4\;,\label{2.9}
\end{align}
where 
\begin{align*}
&J_1=c\,\|\Delta_j u^{(n+1)}\|_{L^2}\sum_{m\le j-1}2^{(1+\frac{d}{2})m}\|\Delta_m u^{(n)}\|_{L^2}\;,\\&J_2=c\,\|\Delta_j u^{(n)}\|_{L^2}\sum_{m\le j}2^{(1+\frac{d}{2})m}\|\Delta_m u^{(n+1)}\|_{L^2}\;,\\&J_3=c\,2^j\sum_{k\ge j-1} 2^{\frac{d}{2}k}\|\widetilde{\Delta}_k u^{(n+1)}\|_{L^2}\|\Delta_k u^{(n)}\|_{L^2}\;,\\& J_4=c\, 2^j \|\Delta_j w^{(n)}\|_{L^2}\;.
\end{align*}
Integrating  (\ref{2.9}) in time yields
\begin{eqnarray}
\|\Delta_j u^{(n+1)}\|_{L^2}\le e^{-c_02^{2\alpha j}t}\|\Delta_j u_0^{(n+1)}\|_{L^2}+\int_0^t e^{-c_0 2^{2\alpha j}(t-\tau)}(J_1+\cdots+J_4)\,d\tau~.\label{33.5}
\end{eqnarray}
Multiplying (\ref{33.5}) by $2^{(1+\frac{d}{2}-2\alpha)j}$ and summing over $j,$ we obtain
\begin{align}
\|u^{(n+1)}(t)\|_{B_{2,1}^{1+\frac{d}{2}-2\alpha}}\le \|u_0^{(n+1)}\|_{B_{2,1}^{1+\frac{d}{2}-2\alpha}}+\sum_{j\ge -1} 2^{(1+\frac{d}{2}-2\alpha)}\int_0^{t}e^{-c_0 2^{2\alpha j}(t-\tau)}(J_1+\cdots +J_4)\,d\tau\;.\label{3.6}
\end{align}
The terms on the right hand side of (\ref{3.6}) can be estimated as follows using the simple bound
\begin{align*}
e^{-c_0 2^{2\alpha j}(t-\tau)}\le 1\;.
\end{align*}
Recalling the definition of $J_1$ above and using the inductive assumption on $u^{(n)},$ we have for any $t\le T,$
\begin{align*}
\sum_{j\ge -1} 2^{(1+\frac{d}{2}-2\alpha)j}&\int_0^t e^{-c_02^{2\alpha j}(t-\tau)}J_1\,d\tau\\&\le c \int_0^t\sum_{j\ge -1} 2^{(1+\frac{d}{2}-2\alpha)j}|\Delta_j u^{(n+1)}\|_{L^2}\sum_{m\le j-1}2^{(1+\frac{d}{2})m}\|\Delta_m u^{(n)}(\tau)\|_{L^2}\,d\tau
\\&\le c \int_0^t\|u^{(n+1)}(\tau)\|_{B_{2,1}^{1+\frac{d}{2}-2\alpha}}\|u^{(n)}(\tau)\|_{B_{2,1}^{1+\frac{d}{2}}}\,d\tau\\&\le c\,\|u^{(n+1)}\|_{L^{\infty}(0,t,B_{2,1}^{1+\frac{d}{2}-2\alpha})}\|u^{(n)}\|_{L^1(0,t,B_{2,1}^{1+\frac{d}{2}})}\\&\le c\,\|u^{(n+1)}\|_{L^{\infty}(0,T,B_{2,1}^{1+\frac{d}{2}-2\alpha})}\|u^{(n)}\|_{L^1(0,T,B_{2,1}^{1+\frac{d}{2}})}\le c\,\delta\,\|u^{(n+1)}\|_{L^{\infty}(0,T,B_{2,1}^{1+\frac{d}{2}-2\alpha})}~.
\end{align*}
The term involving $J_2$ admits the same bound. In fact, by Young's inequality for series convolution,
\begin{align*}
\sum_{j\ge -1}2^{(1+\frac{d}{2}-2\alpha)j}\int_0^t&e^{-c_02^{2\alpha j}(t-\tau)}J_2\,d\tau\\&\le c\int_0^t \sum_{j\ge -1}2^{(1+\frac{d}{2})j}\|\Delta_j u^{(n)}\|_{L^2}\sum_{m\le j}2^{2\alpha(m-j)}2^{(1+\frac{d}{2}-2\alpha)m}\|\Delta_mu^{(n+1)}(\tau)\|_{L^2}\,d\tau\\&\le c\int_0^t \|u^{(n)}(\tau)\|_{B_{2,1}^{1+\frac{d}{2}}}\|u^{(n+1)}(\tau)\|_{B_{2,1}^{1+\frac{d}{2}-2\alpha}}\,d\tau \\&\le c\|u^{(n)}\|_{L^1(0,T,B_{2,1}^{1+\frac{d}{2}})}\|u^{(n+1)}\|_{L^{\infty}(0,T, B_{2,1}^{1+\frac{d}{2}-2\alpha})}\le c\,\delta\,\|u^{(n+1)}\|_{L^{\infty}(0,T,B_{2,1}^{1+\frac{d}{2}-2\alpha})}~.\end{align*}
The estimate for the term with $J_3$ is also similar, 
\begin{align*}
\sum_{j\ge -1} 2^{(1+\frac{d}{2}-2\alpha)j}\int_0^t &e^{-c_02^{2\alpha j}(t-\tau)}J_3\,d\tau\\&=\int_0^t \sum_{j\ge -1}2^{(1+\frac{d}{2}-2\alpha)j}\sum_{k\ge j-1}c\,2^{j}2^{\frac{d}{2}k}\|\widetilde{\Delta}_ku^{(n+1)}\|_{L^2}\|\Delta_k u^{(n)}\|_{L^2}\,d\tau\\&=c\int_0^t \sum_{j\ge-1}\sum_{k\ge j-1}2^{(2+\frac{d}{2}-2\alpha)(j-k)}2^{(1+\frac{d}{2})k}\|\Delta_k u^{(n)}\|_{L^2}2^{(1+\frac{d}{2}-2\alpha)k}\|\widetilde{\Delta}_ku^{(n+1)}\|_{L^2}\,d\tau\\&\le c\int_0^t \|u^{(n)}(\tau)\|_{B_{2,1}^{1+\frac{d}{2}}}\|u^{(n+1)}(\tau)\|_{B_{2,1}^{1+\frac{d}{2}-2\alpha}}\,d\tau\\&\le c\,\|u^{(n)}\|_{L^1(0,T, B_{2,1}^{1+\frac{d}{2}})}\|u^{(n+1)}\|_{L^{\infty}(0,T, B_{2,1}^{1+\frac{d}{2}-2\alpha})}\le c\,\,\delta\,\|u^{(n+1)}\|_{L^{\infty}(0,T, B_{2,1}^{1+\frac{d}{2}-2\alpha})}~.\end{align*}
It remains to bound the term with $J_4,$ 
\begin{align*}
\sum_{j\ge -1}2^{(1+\frac{d}{2}-2\alpha)j}\int_0^t& e^{-c_02^{2\alpha j}(t-\tau)}J_4\,d\tau\\&=\sum_{j\ge -1}2^{(1+\frac{d}{2}-2\alpha)j}\int_0^t e^{-c_02^{2\alpha j}(t-\tau)} c\,2^j\|\Delta_j w^{(n)}\|_{L^2}\,d\tau\\&\le \sum_{j\ge -1}2^{(1+\frac{d}{2}-2\alpha)j}\int_0^t  c\,2^j\,\|\Delta_j w^{(n)}\|_{L^2}\,d\tau= c\int_0^t\sum_{j\ge -1}2^{(2+\frac{d}{2}-2\alpha)j}\|\Delta_j w^{(n)}\|_{L^2}\,d\tau\\&\underset{\text{since}\;\alpha\ge \frac{1}{2}}{\le} c\int_0^t\sum_{j\ge -1}2^{(1+\frac{d}{2})j}\|\Delta_jw^{(n)}\|_{L^2}\,d\tau =c\,\|w^{(n)}\|_{L^1(0,t,B_{2,1}^{1+\frac{d}{2}})}\le c\,\delta~.
\end{align*}
Collecting the bounds above and inserting them in (\ref{3.6}), we find for any $t\le T$ 
\begin{align*}
\|u^{(n+1)}(t)\|_{B_{2,1}^{1+\frac{d}{2}-2\alpha}}\le \|u_0^{(n+1)}\|_{B_{2,1}^{1+\frac{d}{2}-2\alpha}}+c\,\delta\,\|u^{(n+1)}\|_{L^{\infty}(0,T,B_{2,1}^{1+\frac{d}{2}-2\alpha})}+c\,\delta~.
\end{align*}
Therefore
\begin{align*}
\|u^{(n+1)}(t)\|_{L^{\infty}(0,T,B_{2,1}^{1+\frac{d}{2}-2\alpha})}\le \|u_0^{(n+1)}\|_{B_{2,1}^{1+\frac{d}{2}-2\alpha}}+c\,\delta\,\|u^{(n+1)}\|_{L^{\infty}(0,T,B_{2,1}^{1+\frac{d}{2}-2\alpha})}+c\,\delta~.
\end{align*}
Choosing $\delta$ such that $c\,\delta \le \min(\frac{1}{4},\frac{M}{4})$ we get
\begin{align*}
\|u^{(n+1)}(t)\|_{L^{\infty}(0,T,B_{2,1}^{1+\frac{d}{2}-2\alpha})}\le \frac{M}{2}+\frac{1}{4}\|u^{(n+1)}\|_{L^{\infty}(0,T,B_{2,1}^{1+\frac{d}{2}-2\alpha})}+\frac{M}{4}~,
\end{align*}
which implies
\begin{align*}
\|u^{(n+1)}(t)\|_{L^{\infty}(0,T,B_{2,1}^{1+\frac{d}{2}-2\alpha})}\le M~.
\end{align*}
\subsubsection{The estimate of $\|u^{(n+1)}\|_{L^1(0,T,B^{1+\frac{d}{2}})}$.} We multiply (\ref{33.5}) by $2^{(1+\frac{d}{2})j}$, sum over $j$ and integrate in time to obtain
\begin{align}
\|u^{(n+1)}\|_{L^1(0,T,B_{2,1}^{1+\frac{d}{2}})}\le \int_0^T\sum_{j\ge -1}&2^{(1+\frac{d}{2})j}e^{-c_02^{2\alpha j}t}\|\Delta_j u_0^{(n+1)}\|_{L^2}\,dt\nonumber\\&+\int_0^T\int_0^s\sum_{j\ge -1}2^{(1+\frac{d}{2})j}e^{-c_02^{2\alpha j}(s-\tau)}(J_1+\cdots +J_4)\,d\tau\,ds\;.
\label{526}\end{align}
We estimate the terms on the right hand side of (\ref{526}) and start with the first term.
\begin{align*}
\int_0^T\sum_{j\ge -1}2^{(1+\frac{d}{2})j}e^{-c_02^{2\alpha j}t}\|\Delta_ju_0^{n+1}\|_{L^2}\,dt=c\sum_{j\ge -1}2^{(1+\frac{d}{2}-2\alpha)j}(1-e^{-c_02^{2\alpha j}T})\|\Delta_j u_0^{(n+1)}\|_{L^2}\;.
\end{align*}
Since $u_0\in B_{2,1}^{(1+\frac{d}{2}-2\alpha)},$ then by Dominated convergence Theorem
\begin{align*}
\lim_{T\to 0}\sum_{j\ge -1}2^{(1+\frac{d}{2}-2\alpha)j}(1-e^{-c_02^{2\alpha j}T})\|\Delta_j u_0^{(n+1)}\|_{L^2}=0\;.
\end{align*}
Therefore, we can choose $T$ sufficiently small such that
\begin{align*}
\int_0^T\sum_{j\ge -1}2^{(1+\frac{d}{2})j}e^{-c_02^{2\alpha j}t}\|\Delta_j u_0^{(n+1)}\|_{L^2}\,dt\le \frac{\delta}{4}\;.
\end{align*}
Applying Young's inequality for the time convolution, we have 
\begin{align*}
\int_0^T&\sum_{j\ge -1}2^{(1+\frac{d}{2})j}\int_0^se^{-c_02^{2\alpha j}(s-\tau)} J_1\,d\tau\,ds\\&=c\int_0^T\sum_{j\ge -1}2^{(1+\frac{d}{2})j}\int_0^se^{-c_02^{2\alpha j}(s-\tau)} \|\Delta_j u^{(n+1)}(\tau)\|_{L^2}\sum_{m\le j -1}2^{(1+\frac{d}{2})m} \|\Delta_m u^{(m)}\|_{L^2}\,d\tau\,ds\\&\le c\sum_{j\ge -1}2^{(1+\frac{d}{2})j}\int_0^T\|\Delta_j u^{(n+1)}(\tau)\|_{L^2}\sum_{m\le j -1}2^{(1+\frac{d}{2})m} \|\Delta_m u^{(m)}\|_{L^2}\,d\tau \int_0^Te^{-c_02^{2\alpha j}s}\,ds\end{align*}
Using the fact that there exists $c_2> 0$ satisfying for $j\ge 0\;,$
\begin{align}
\int_0^se^{-c_02^{2\alpha j}s}ds\le c\,2^{-2\alpha j}(1-e^{-c_2T})\;,\label{sou}
\end{align}
we get
\begin{align*}\int_0^T\sum_{j\ge -1}2^{(1+\frac{d}{2})j}&\int_0^se^{-c_02^{2\alpha j}(s-\tau)} J_1\,d\tau\,ds\\&\le c\,(1-e^{-c_2T})\int_0^T\sum_j2^{(1+\frac{d}{2}-2\alpha)j}\|\Delta_j u^{(n+1)}(\tau)\|_{L^2}\sum_{m\le j -1}2^{(1+\frac{d}{2})m} \|\Delta_m u^{(m)}\|_{L^2}\,d\tau\\&\le c\,(1-e^{-c_2T})\|u^{(n+1)}\|_{L^{\infty}(0,T,B_{2,1}^{1+\frac{d}{2}-2\alpha})}\underbrace{\|u^{(n)}\|_{L^{1}(0,T,B_{2,1}^{1+\frac{d}{2}})}}_{\le \delta}\\&\le c\,\delta\,(1-e^{-c_2T})\|u^{(n+1)}\|_{L^{\infty}(0,T,B_{2,1}^{1+\frac{d}{2}-2\alpha})}\;.
\end{align*}
The terms with $J_2$ and $J_3$ can be similarly estimated and obey the same bound.
\begin{align*}\int_0^T&\sum\limits_{j\ge -1}2^{(1+\frac{d}{2})j}\int_0^se^{-c_02^{2\alpha j}(s-\tau)} J_2\,d\tau\,ds
\\&=c\int_0^T\sum\limits_{j\ge -1}2^{(1+\frac{d}{2})j}\int_0^se^{-c_02^{2\alpha j}(s-\tau)} \|\Delta_ju^{(n)}\|_{L^2}\sum\limits_{m\le j}2^{(1+\frac{d}{2})m}\|\Delta_mu^{(n+1)}\|_{L^2}\,d\tau\,ds
\\&\le c\sum\limits_{j\ge -1}2^{(1+\frac{d}{2})j}\int_0^T\|\Delta_ju^{(n)}\|_{L^2}\sum\limits_{m\le j}2^{(1+\frac{d}{2})m}\|\Delta_mu^{(n+1)}\|_{L^2}\,d\tau\int_0^Te^{-c_02^{2\alpha j}s}\,ds\;.\end{align*}
Owing to (\ref{sou}) and the above inequality,
\begin{align*}\int_0^T\sum\limits_{j\ge -1}2^{(1+\frac{d}{2})j}&\int_0^se^{-c_02^{2\alpha j}(s-\tau)} J_2\,d\tau\,ds
\\&\le c (1-e^{-c_2T})\int_0^T\sum\limits_{j\ge -1}2^{(1+\frac{d}{2}-2\alpha)j}\|\Delta_ju^{(n)}\|_{L^2}\sum\limits_{m\le j-1}2^{(1+\frac{d}{2})m}\|\Delta_mu^{(n+1)}\|_{L^2}\,d\tau\\&\le c (1-e^{-c_2T})\int_0^T\sum\limits_{j\ge -1}2^{(1+\frac{d}{2}-2\alpha)j}\|\Delta_ju^{(n+1)}\|_{L^2}\sum\limits_{j\ge -1}2^{(1+\frac{d}{2})j}\|\Delta_ju^{(n)}\|_{L^2}\,d\tau \\&\le c(1-e^{-c_2T})\|u^{(n+1)}\|_{L^{\infty}(0,T,B_{2,1}^{1+\frac{d}{2}-2\alpha})}\underbrace{\|u^{(n)}\|_{L^{1}(0,T,B_{2,1}^{1+\frac{d}{2}})}}_{\le\delta}\\&\le c(1-e^{-c_2T})\,\delta\,\|u^{(n+1)}\|_{L^{\infty}(0,T,B_{2,1}^{1+\frac{d}{2}-2\alpha})}\;.\end{align*}
Similarly
\begin{align*}\int_0^T\sum\limits_{j\ge -1}2^{(1+\frac{d}{2})j}&\int_0^se^{-c_02^{2\alpha j}(s-\tau)} J_3\,d\tau\,ds
\\&=c\int_0^T\sum\limits_{j\ge -1}2^{(2+\frac{d}{2})j}\int_0^se^{-c_02^{2\alpha j}(s-\tau)}\sum\limits_{k\ge j-1}2^{\frac{d}{2}k}\|\Delta_ku^{(n)}\|_{L^2}\|\widetilde{\Delta}_ku^{(n+1)}\|_{L^2}\,d\tau\,ds\\&\le c\sum\limits_{j\ge -1}2^{(2+\frac{d}{2})j}\int_0^T\sum\limits_{k\ge j-1}2^{\frac{d}{2}k}\|\Delta_ku^{(n)}\|_{L^2} \|\widetilde{\Delta}_ku^{(n+1)}\|_{L^2}d\tau\int_0^Te^{-c_02^{2\alpha j}s}\,ds~.
\end{align*}
Then, due to (\ref{sou}),
\begin{align*}\int_0^T\sum\limits_{j\ge -1}2^{(1+\frac{d}{2})j}&\int_0^se^{-c_02^{2\alpha j}(s-\tau)} J_3\,d\tau\,ds
\\&\le c\,(1-e^{-c_2T})\int_0^T\sum\limits_{j\ge -1}2^{(2+\frac{d}{2}-2\alpha)j}\sum\limits_{k\ge j-1}2^{\frac{d}{2}k}\|\Delta_ku^{(n)}\|_{L^2}\|\widetilde{\Delta}_ku^{(n+1)}\|_{L^2}\,d\tau
\\&\le c\,(1-e^{-c_2T})\int_0^T\sum\limits_{j\ge -1}2^{(2+\frac{d}{2}-2\alpha)j}\|\widetilde{\Delta}_ju^{(n+1)}\|_{L^2}\sum\limits_{j\ge -1}2^{(\frac{d}{2})j}\|\Delta_ju^{(n)}\|_{L^2}\,d\tau
\\&\le c\,(1-e^{-c_2T})\int_0^T\sum\limits_{j\ge -1}2^{(1+\frac{d}{2}-2\alpha)j}\|\widetilde{\Delta}_ju^{(n+1)}\|_{L^2}\sum\limits_{j\ge -1}2^{(1+\frac{d}{2})j}\|\Delta_ju^{(n)}\|_{L^2}\,d\tau 
\\&\le c\,(1-e^{-c_2T})\|u^{(n+1)}\|_{L^{\infty}(0,T,B_{2,1}^{1+\frac{d}{2}-2\alpha})}\underbrace{\|u^{(n)}\|_{L^{1}(0,T,B_{2,1}^{1+\frac{d}{2}})}}_{\le\delta}\\&\le c\,(1-e^{-c_2T})\,\delta\,\|u^{(n+1)}\|_{L^{\infty}(0,T,B_{2,1}^{1+\frac{d}{2}-2\alpha})}\;.\end{align*}
Now, for the term with $J_4$ we write
\begin{align*}\int_0^T\sum\limits_{j\ge -1}2^{(1+\frac{d}{2})j}\int_0^se^{-c_02^{2\alpha j}(s-\tau)} J_4\,d\tau\,ds
&=c\int_0^T\sum\limits_{j\ge -1}2^{(2+\frac{d}{2})j}\int_0^se^{-c_02^{2\alpha j}(s-\tau)} \|\Delta_jw^{(n)}\|_{L^2}\,d\tau\,ds
\\&\le c\sum\limits_{j\ge -1}2^{(2+\frac{d}{2})j}\int_0^T\|\Delta_jw^{(n)}\|_{L^2}\,d\tau\int_0^Te^{-c_02^{2\alpha j}s}\,ds
\\&\le c\,(1-e^{-c_2T})\int_0^T\sum\limits_{j\ge -1}2^{(2+\frac{d}{2}-2\alpha)}\|\Delta_jw^{(n)}\|_{L^2}\,d\tau
\\&\underset{\text{since}\,\alpha\ge\frac{1}{2}}{\le} c\,(1-e^{-c_2T})\underbrace{\int_0^T\sum\limits_{j\ge -1}2^{(1+\frac{d}{2})j}\|\Delta_jw^{(n)}\|_{L^2}\,d\tau
}_{=\|w^{(n)}\|_{L^{1}(0,T,B_{2,1}^{1+\frac{d}{2}})}\le \delta}\\&\le c\,(1-e^{-c_2T})\delta\;.\end{align*}
Collecting the estimates above leads to
\begin{align*}
\|u^{(n+1)}\|_{L^{1}(0,T,B_{2,1}^{1+\frac{d}{2}})}&\le \frac{\delta}{4}+c\,\delta(1-e^{-c_2T})\|u^{(n+1)}\|_{L^{\infty}(0,T,B_{2,1}^{1+\frac{d}{2}})}+c\,(1-e^{-c_2T})\delta\\&\le \frac{\delta}{4}+c\,\delta(1-e^{-c_2T})M+c\,(1-e^{-c_2T})\delta\;.
\end{align*}
Choosing T sufficiently small such that 
$c\,(1-e^{-c_2T})\le \min(\frac{1}{4M} , \frac{1}{2})$ we get
\begin{align*}
\|u^{(n+1)}\|_{L^{1}(0,T,B_{2,1}^{1+\frac{d}{2}})}\le \frac{\delta}{4}+\frac{\delta}{4}+\frac{\delta}{2}=\delta\;.
\end{align*}
\subsubsection{The estimate of $w^{(n+1)}$ in $L^\infty(0, T, B_{2,1}^{1+\frac{d}{2}-2\beta}(\mathbb{R}^d)).$} We apply $\Delta_j$ to the third equation in (\ref{2.6}) and then dotting with $\Delta_jw^{(n+1)},$ we obtain
\begin{align}
\frac{1}{2}\frac{d}{dt}\|\Delta_jw^{(n+1)}\|_{L^2}^2+(c_1&2^{2\beta j}+4k)\|\Delta_jw^{(n+1)}\|_{L^2}^2\nonumber\\&\le -2k\int \Delta_j(\nabla\times u^{(n)})\Delta_jw^{(n+1)}\,dx\nonumber\\&-\int \Delta_j( u^{(n)}\cdot\nabla w^{(n+1)})\Delta_jw^{(n+1)}\,dx\nonumber\\&=B_1+B_2\;,\label{520}
\end{align}
where 
\begin{align*}
&B_1=-2k\int \Delta_j(\nabla\times u^{(n)})\Delta_jw^{(n+1)}\,dx\;,\\
&B_2=-\int \Delta_j( u^{(n)}\cdot\nabla w^{(n+1)})\Delta_jw^{(n+1)}\,dx\;.
\end{align*}
By H\"older's inequality and Bernstein's inequality,
\begin{align*}
|B_1|&=|-2k\int \Delta_j(\nabla\times u^{(n)})\Delta_jw^{(n+1)}\,dx|\\&\le 2k\,\|\Delta_j(\nabla\times u^{(n)})\|_{L^2}\|\Delta_jw^{(n+1)}\|_{L^2}\\&\le c\,2^j\,\|\Delta_j u^{(n)}\|_{L^2}\|\Delta_jw^{(n+1)}\|_{L^2}\;.
\end{align*}
By Lemma \ref{lem2.3},
\begin{align*}
|B_2|&=|-\int \Delta_j( u^{(n)}\cdot\nabla w^{(n+1)})\Delta_jw^{(n+1)}\,dx|\\&\le c\,\|\Delta_j w^{(n+1)}\|_{L^2}^2\sum\limits_{m\le j-1}2^{(1+\frac{d}{2})m}\|\Delta_mu^{(n)}\|_{L^2}\\
& +c\,\|\Delta_j w^{(n+1)}\|_{L^2}\|\Delta_j u^{(n)}\|_{L^2}\sum\limits_{m\le j}2^{(1+\frac{d}{2})m}\|\Delta_mw^{(n+1)}\|_{L^2}\\&+c\,\|\Delta_j w^{(n+1)}\|_{L^2}2^j\sum\limits_{k\ge j-1}2^{\frac{d}{2}k}\|\widetilde{\Delta}_k w^{(n+1)}\|_{L^2}\|\Delta_ku^{(n)}\|_{L^2}\;.
\end{align*}
Inserting the estimates above in (\ref{520}) and eliminating $\|\Delta_j w^{(n+1)}\|_{L^2}$ from both sides of the inequality, we obtain 
\begin{align}
\frac{d}{dt}\|\Delta_jw^{(n+1)}\|_{L^2}+(c_12^{2\beta j}+4k)\|\Delta_jw^{(n+1)}\|_{L^2}\le K_1+K_2+K_3+K_4\;,
\label{521}\end{align}
where 
\begin{align*}
&K_1=c\,2^{j}\,\|\Delta_j u^{(n)})\|_{L^2}\;,\\
&K_2=c\,\|\Delta_j w^{(n+1)}\|_{L^2}\sum\limits_{m\le j-1}2^{(1+\frac{d}{2})m}\|\Delta_mu^{(n)}\|_{L^2}\;,\\
&K_3=c\,\|\Delta_j u^{(n)}\|_{L^2}\sum\limits_{m\le j}2^{(1+\frac{d}{2})m}\|\Delta_mw^{(n+1)}\|_{L^2}\;,\\
&K_4=c\,\sum\limits_{k\ge j-1}2^j2^{\frac{d}{2}k}\|\widetilde{\Delta}_k w^{(n+1)}\|_{L^2}\|\Delta_ku^{(n)}\|_{L^2}\;.
\end{align*}
Integrating (\ref{521}) in time yields, for any $t\le T,$
\begin{eqnarray}
\|\Delta_jw^{(n+1)}\|_{L^2}\le e^{-(c_12^{2\beta j})t}\|\Delta_jw_0^{(n+1)}\|_{L^2}+\int_0^te^{-(c_12^{2\beta j})(t-\tau)}(K_1+\cdots+K_4)\,d\tau\;.\label{3.12}
\end{eqnarray}
\\
Multiplying (\ref{3.12}) by $2^{(1+\frac{d}{2}-2\beta) j}$ and summing over $j$, we have
\begin{align}
\|w^{(n+1)}\|_{B_{2,1}^{1+\frac{d}{2}-2\beta}}\le \|w^{(n+1)}_0\|_{B_{2,1}^{1+\frac{d}{2}-2\beta}}+\sum_{j\ge -1}\int_0^te^{-(c_12^{2\beta j})(t-\tau)}2^{(1+\frac{d}{2}-2\beta)j}(K_1+\cdots+K_4)\,d\tau\;.\label{3.13}
\end{align}
The terms containing $K_1$ through $K_4$ on the right hand side of $(\ref{3.13})$ can be bounded suitably as follows. We start with the term with $K_1$,
\begin{align*}
\sum_{j\ge -1}\int_0^t2^{(1+\frac{d}{2}-2\beta)j}K_1\,d\tau&=\int_0^t\sum_{j\ge -1}c\,2^{(2+\frac{d}{2}-2\beta)j}\|\Delta_ju^{(n)}(\tau)\|_{L^2}\,d\tau\\&\underset{\text{since}\;\beta\ge\frac{1}{2}}{\le} c\underbrace{\int_0^t\sum_{j\ge -1}\,2^{(1+\frac{d}{2})j}\|\Delta_ju^{(n)}(\tau)\|_{L^2}\,d\tau}_{ \|u^{(n)}\|_{L^{1}(0,T,B_{2,1}^{1+\frac{d}{2}})}\le \delta}\le c\,\delta\,.
\end{align*}
Similarly the term with $K_2$ is bounded by 
\begin{align*}
\sum_{j\ge -1}2^{(1+\frac{d}{2}-2\beta)j}\int_0^tK_2\,d\tau&=c\int_0^t\sum_{j\ge -1}2^{(1+\frac{d}{2}-2\beta)j}\|\Delta_jw^{(n+1)}\|_{L^2} \sum_{m\le j-1}2^{(1+\frac{d}{2})m}\|\Delta_mu^{(n)}\|_{L^2} \,d\tau\\&\le c\,\|w^{(n+1)}\|_{L^{\infty}(0,T,B_{2,1}^{1+\frac{d}{2}-2\beta})}\|u^{(n)}\|_{L^{1}(0,T,B_{2,1}^{1+\frac{d}{2}})}\\&\le c\,\delta\,\|w^{(n+1)}\|_{L^{\infty}(0,T,B_{2,1}^{1+\frac{d}{2}-2\beta})}\;.
\end{align*}
The terms related to $K_3$  and $K_4$ obey also the same bound,
\begin{align*}
\sum_{j\ge -1}2^{(1+\frac{d}{2}-2\beta)j}\int_0^tK_3\,d\tau&=c\int_0^t\sum_{j\ge -1}2^{(1+\frac{d}{2}-2\beta)j}\|\Delta_ju^{(n)}\|_{L^2} \sum_{m\le j}2^{(1+\frac{d}{2})m}\|\Delta_mw^{(n+1)}\|_{L^2}\,d\tau \\&\le c\,\|w^{(n+1)}\|_{L^{\infty}(0,T,B_{2,1}^{1+\frac{d}{2}-2\beta})}\|u^{(n)}\|_{L^{1}(0,T,B_{2,1}^{1+\frac{d}{2}})}\\&\le c\,\delta\,\|w^{(n+1)}\|_{L^{\infty}(0,T,B_{2,1}^{1+\frac{d}{2}-2\beta})}\;.
\end{align*}
For the term with $K_4$ we write
\begin{align*}
\sum_{j\ge -1}2^{(1+\frac{d}{2}-2\beta)j}\int_0^tK_4\,d\tau&=c\int_0^t\sum_{j\ge -1}\sum_{k\ge j-1}2^{(2+\frac{d}{2}-2\beta)j}2^{\frac{d}{2}k}\|\widetilde{\Delta}_kw^{(n+1)}\|_{L^2} \|\Delta_ku^{(n)}\|_{L^2}\,d\tau\\&\le c\int_0^t\sum_{j\ge -1}\sum_{j\ge -1}2^{(2+\frac{d}{2}-2\beta)j}2^{\frac{d}{2}j}\|\widetilde{\Delta}_jw^{(n+1)}\|_{L^2} \|\Delta_ju^{(n)}\|_{L^2}\,d\tau\\&=c\int_0^t\sum_{j\ge -1}\sum_{j\ge -1}2^{(1+\frac{d}{2}-2\beta)j}2^{(1+\frac{d}{2})j}\|\widetilde{\Delta}_jw^{(n+1)}\|_{L^2} \|\Delta_ju^{(n)}\|_{L^2}\,d\tau\\&\le c\,\|w^{(n+1)}\|_{L^{\infty}(0,T,B_{2,1}^{1+\frac{d}{2}-2\beta})}\|u^{(n)}\|_{L^{1}(0,T,B_{2,1}^{1+\frac{d}{2}})}\\&\le c\,\delta\|w^{(n+1)}\|_{L^{\infty}(0,T,B_{2,1}^{1+\frac{d}{2}-2\beta})}\;.
\end{align*}
Collecting the estimates and inserting them in (\ref{3.13}), we obtain for any $t\le T$
\begin{align*}
\|w^{(n+1)}(t)\|_{B_{2,1}^{1+\frac{d}{2}-2\beta}}\le \|&w_0^{(n+1)}\|_{B_{2,1}^{1+\frac{d}{2}-2\beta}}+c\,\delta+c\,\delta \|w^{(n+1)}\|_{L^{\infty}(0,T,B_{2,1}^{1+\frac{d}{2}-2\beta})}~.
\end{align*}
Choosing $c\,\delta\le\min(\frac{1}{4},\frac{M}{4}),$ we get
\begin{align*}
\|w^{(n+1)}(t)\|_{L^{\infty}(0,T,B_{2,1}^{1+\frac{d}{2}-2\beta})}\le\frac{M}{2}+\frac{M}{4}+\frac{1}{4}\|w^{(n+1)}(t)\|_{L^{\infty}(0,T,B_{2,1}^{1+\frac{d}{2}-2\beta})}~,
\end{align*}
which implies
\begin{align*}
\|w^{(n+1)}(t)\|_{L^{\infty}(0,T,B_{2,1}^{1+\frac{d}{2}-2\beta})}\le M~.
\end{align*}
\subsubsection{The estimate of $\|w^{(n+1)}(t)\|_{L^{1}(0,T,B_{2,1}^{1+\frac{d}{2}})}.$} We recall (\ref{3.12})
\begin{align*}
\|\Delta_jw^{(n+1)}\|_{L^2}\le e^{-(c_12^{2\beta j})t}\|\Delta_jw^{(n+1)}_0\|_{L^2}+\int_0^te^{-c_12^{2\beta j}(t-\tau)}(K_1+\cdots+K_4)\,d\tau~.
\end{align*}
We multiply by $2^{(1+\frac{d}{2})j},$ sum over $j$ and integrate in time to get 
\begin{align}
\|w^{(n+1)}\|_{L^1(0,T,B_{2,1}^{1+\frac{d}{2}})}&\le\int_0^T\sum_{j\ge -1}2^{(1+\frac{d}{2})j}e^{-c_12^{2\beta j}t} \|\Delta_jw^{(n+1)}_0\|_{L^2}\,dt\nonumber\\&+ \int_0^T\sum_{j\ge -1}2^{(1+\frac{d}{2})j}\int_0^se^{-c_12^{2\beta j}(s-\tau)}(K_1+\cdots +K_4)\,d\tau\,ds~.\label{522}
\end{align}
Clearly
\begin{align*}
\int_0^T\sum_{j\ge -1}2^{(1+\frac{d}{2})j}&e^{-c_12^{2\beta j}t} \|\Delta_jw^{(n+1)}_0\|_{L^2}\,dt=c\sum_{j\ge -1}2^{(1+\frac{d}{2}-2\beta)j}(1-e^{-c_12^{2\beta j}T})\|\Delta_jw^{(n+1)}_0\|_{L^2}\;.
\end{align*}
Since $w_0\in B_{2,1}^{1+\frac{d}{2}-2\beta},$ we have by the Dominated Convergence Theorem,
\begin{align*}
\lim_{T\to 0}\sum_{j\ge -1}2^{(1+\frac{d}{2}-2\beta)j}(1-e^{-c_12^{2\beta j}T})\|\Delta_jw^{(n+1)}_0\|_{L^2}=0\;.
\end{align*}
Therefore, we can choose $T$ sufficiently small 
\begin{align*}
\int_0^T\sum_{j\ge -1}2^{(1+\frac{d}{2})j}e^{-c_12^{2\beta j}t}\|\Delta_jw^{(n+1)}_0\|_{L^2}\,dt\le \frac{\delta}{2}\;.
\end{align*}
Applying Young's inequality for the time convolution, the term with $K_1$ is bounded by 
\begin{align*}
\int_0^T\sum_{j\ge -1}2^{(1+\frac{d}{2})j}\int_0^se^{-c_12^{2\beta j}(s-\tau)} K_1\,d\tau \,ds&\le c\int_0^T\sum_{j\ge -1}2^{(2+\frac{d}{2})j}\int_0^se^{-c_12^{2\beta j}(s-\tau)} \|\Delta_ju^{(n)}\|_{L^2}\,d\tau\,ds\\&\le c\,\sum_{j\ge -1}2^{(2+\frac{d}{2})j}\int_0^T \|\Delta_ju^{(n)}\|_{L^2}\,d\tau\cdot \int_0^Te^{-c_12^{2\beta j}s}\,ds.
\end{align*}
Using the fact that there exists $ c_3>0$ satisfying for all $j\ge 0\;,$
\begin{align}
\int_0^Te^{-c_12^{2\beta j}s} ds\le c\, 2^{-2\beta j}(1-e^{-c_3T})\;,\label{sou1}
\end{align}
we get
\begin{align*}
\int_0^T\sum_{j\ge -1}2^{(1+\frac{d}{2})j}\int_0^se^{-c_12^{2\beta j}(s-\tau)} K_1\,d\tau \,ds&\le c\,(1-e^{-c_3T})\int_0^T\sum_{j\ge -1}2^{(2+\frac{d}{2}-2\beta)j}\|\Delta_ju^{(n)}(\tau)\|_{L^2}\,d\tau\\&\le c\,(1-e^{-c_3T})\underbrace{\int_0^T\sum_{j\ge -1}2^{(1+\frac{d}{2})j}\|\Delta_ju^{(n)}(\tau)\|_{L^2}\,d\tau}_{=\|u^{(n+1)}\|_{L^{1}(0,T,B_{2,1}^{1+\frac{d}{2}})}\le \delta}\\&\le c\,(1-e^{-c_3T})\delta\;.
\end{align*}
Similarly by applying Young's inequality for the time convolution, the term with $K_2$ is bounded by
\begin{align*}
\int_0^T\sum_{j\ge -1}&2^{(1+\frac{d}{2})j}\int_0^se^{-c_12^{2\beta j}(s-\tau)} K_2\,d\tau\,ds\\&=c\,\int_0^T\sum_{j\ge -1}2^{(1+\frac{d}{2})j}\int_0^se^{-c_12^{2\beta j}(s-\tau)} c\,\|\Delta_jw^{(n+1)}\|_{L^2}\sum_{m\le j-1}2^{(1+\frac{d}{2})m}\|\Delta_mu^{(n)}(\tau)\|_{L^2}\,d\tau\,ds\\&\le c\sum_{j\ge -1}2^{(1+\frac{d}{2})j}\int_0^T \|\Delta_jw^{(n+1)}\|_{L^2}\sum_{m\le j-1}2^{(1+\frac{d}{2})m}\|\Delta_mu^{(n)}(\tau)\|_{L^2}\,d\tau\Big(\int_0^Te^{-c_12^{2\beta j}s} ds\Big)\\&\le c\,(1-e^{-c_3T})\int_0^T\sum_{j\ge -1}2^{(1+\frac{d}{2}-2\beta)j}\|\Delta_jw^{(n+1)}(\tau)\|_{L^2}\sum_{m\le j-1}2^{(1+\frac{d}{2})m}\|\Delta_mu^{(n)}(\tau)\|_{L^2}\,d\tau\\&\le c\,(1-e^{-c_3T})\underbrace{\|w^{(n+1)}\|_{L^{\infty}(0,T,B_{2,1}^{1+\frac{d}{2}-2\beta})}}_{\le M}\underbrace{\|u^{(n)}\|_{L^{1}(0,T,B_{2,1}^{1+\frac{d}{2}})}}_{\le \delta}\\&\le c\,(1-e^{-c_3T})\delta M\;.
\end{align*}
The terms involving  $K_3$ and $K_4$ obey also the same bound,
\begin{align*}
\int_0^T\sum_{j\ge -1}&2^{(1+\frac{d}{2})j}\int_0^se^{-c_12^{2\beta j}(s-\tau)} K_3\,d\tau\,ds\\&=c\int_0^T\sum_{j\ge -1}2^{(1+\frac{d}{2})j}\int_0^se^{-c_12^{2\beta j}(s-\tau)} c\,\|\Delta_ju^{(n)}(\tau)\|_{L^2}\Big(\sum_{m\le j}2^{(1+\frac{d}{2})m}\|\Delta_mw^{(n+1)}(\tau)\|_{L^2}\Big)\,d\tau ds\\&\le c\sum_{j\ge -1}2^{(1+\frac{d}{2})j}\int_0^T \|\Delta_ju^{(n)}(\tau)\|_{L^2}\sum_{m\le j}2^{(1+\frac{d}{2})m}\|\Delta_mw^{(n+1)}(\tau)\|_{L^2}\,d\tau\Big(\int_0^Te^{-c_12^{2\beta j}s} ds\Big)
\end{align*}
Then, owing to (\ref{sou1})
\begin{align*}
\int_0^T\sum_{j\ge -1}&2^{(1+\frac{d}{2})j}\int_0^se^{-c_12^{2\beta j}(s-\tau)} K_3\,d\tau\,ds
\\&\le c\,(1-e^{-c_3T})\int_0^T\sum_{j\ge -1}2^{(1+\frac{d}{2}-2\beta)j}\|\Delta_ju^{(n)}(\tau)\|_{L^2}\sum_{m\le j}2^{(1+\frac{d}{2})m}\|\Delta_mw^{(n+1)}(\tau)\|_{L^2}\,d\tau\\&\le c\,(1-e^{-c_3T})       \int_0^T\sum_{j\ge -1}2^{(1+\frac{d}{2}-2\beta)j}\|\Delta_ju^{(n)}(\tau)\|_{L^2} \sum_{j\ge -1}2^{(1+\frac{d}{2})j}\|\Delta_jw^{(n+1)}(\tau)\|_{L^2}\,d\tau\\&= c\,(1-e^{-c_3T})       \int_0^T\sum_{j\ge -1}2^{(1+\frac{d}{2})j}\|\Delta_ju^{(n)}(\tau)\|_{L^2} \sum_{j\ge -1}2^{(1+\frac{d}{2}-2\beta)j}\|\Delta_jw^{(n+1)}(\tau)\|_{L^2}\,d\tau\\&\le c\,(1-e^{-c_3T})\underbrace{\|w^{(n+1)}\|_{L^{\infty}(0,T,B_{2,1}^{1+\frac{d}{2}-2\beta})}}_{\le M}\underbrace{\|u^{(n)}\|_{L^{1}(0,T,B_{2,1}^{1+\frac{d}{2}})}}_{\le \delta}\\&\le c\,(1-e^{-c_3T})\delta M\;.
\end{align*}
The term containing $K_4$ is bounded by
\begin{align*}
\int_0^T\sum_{j\ge -1}&2^{(1+\frac{d}{2})j}\int_0^se^{-c_12^{2\beta j}(s-\tau)} K_4\,d\tau\,ds\\&=c\int_0^T\sum_{j\ge -1}2^{(1+\frac{d}{2})j}\int_0^se^{-c_12^{2\beta j}(s-\tau)}2^j2^{\frac{d}{2}k}\|\widetilde{\Delta}_kw^{(n+1)}\|_{L^2}\|\Delta_ku^{(n)}\|_{L^2}\,d\tau ds\\&\le c\sum_{j\ge -1}2^{(1+\frac{d}{2})j}\int_0^T\sum_{k\ge j-1}  2^j2^{\frac{d}{2}k}\|\widetilde{\Delta}_kw^{(n+1)}\|_{L^2}\|\Delta_ku^{(n)}\|_{L^2}\,d\tau \Big(\int_0^Te^{-c_12^{2\beta j}s}ds\Big)~.
\end{align*}
Hence, due to (\ref{sou1})
\begin{align*}
\int_0^T\sum_{j\ge -1}&2^{(1+\frac{d}{2})j}\int_0^se^{-c_12^{2\beta j}(s-\tau)} K_4\,d\tau\,ds\\&\le c\,(1-e^{-c_3T})\int_0^T\sum_{j\ge -1}2^{(1+\frac{d}{2}-2\beta)j}2^{j}\sum_{k\ge j-1}2^{\frac{d}{2}k}\|\widetilde{\Delta}_kw^{(n+1)}\|_{L^2}\|\Delta_ku^{(n)}\|_{L^2}\,d\tau \\&= c\,(1-e^{-c_3T})\int_0^T\sum_{j\ge -1}\sum_{k\ge j-1}2^{(1+\frac{d}{2}-2\beta)(j-k)} 2^{(1+\frac{d}{2})k} \|\Delta_ku^{(n)}\|_{L^2}2^{(1+\frac{d}{2}-2\beta)k}  \|\widetilde{\Delta}_kw^{(n+1)}\|_{L^2}\,d\tau\\&= c\,(1-e^{-c_3T})\int_0^T\|u^{(n)}\|_{B_{2,1}^{1+\frac{d}{2}}} \|w^{(n+1)}(\tau)\|_{B_{2,1}^{1+\frac{d}{2}-2\beta}}\,d\tau\\&\le c\,(1-e^{-c_3T})\underbrace{\|u^{(n)}\|_{L^{1}(0,T,B_{2,1}^{1+\frac{d}{2}})}}_{\le \delta}\underbrace{\|w^{(n+1)}\|_{L^{\infty}(0,T,B_{2,1}^{1+\frac{d}{2}-2\beta})}}_{\le M}\\&\le c\,(1-e^{-c_3T})\delta M\;.
\end{align*}
Collecting the estimates above and inserting them in (\ref{522}), we obtain 
\begin{align*}
\|w^{(n+1)}\|_{L^{1}(0,T,B_{2,1}^{1+\frac{d}{2}})}\le \frac{\delta}{2}+c\,(1-e^{-c_3T})\delta +c\,(1-e^{-c_3T})\delta M\;.
\end{align*}
Choosing T sufficiently small such that
$c\,(1-e^{-c_3T})\le\min(\frac{1}{4M},\frac{1}{4})\;,$ we get
\begin{align*}
\|w^{(n+1)}\|_{L^{1}(0,T,B_{2,1}^{1+\frac{d}{2}})}\le \frac{\delta}{2}+\frac{\delta}{4}+\frac{\delta}{4}=\delta\;.
\end{align*}
These uniform bounds allow us to extract a weakly convergent subsequence. That is, there is $(u,w)\in Y$ such that subsequence  of $(u^{n},w^{n})$ (still denoted by $(u^{n},w^{n})$) satisfies 
\begin{align*}
&u^{n}\overset{*}{\rightharpoonup} u \quad\text{in}\quad L^{\infty}(0,T,B_{2,1}^{1+\frac{d}{2}-2\alpha})\;,\\&w^{n}\overset{*}{\rightharpoonup} w \quad\text{in}\quad L^{\infty}(0,T,B_{2,1}^{1+\frac{d}{2}-2\beta})\;.
\end{align*}
In order to show that $(u,w)$ is a weak solution of (\ref{2.1}) we need to further extract a subsequence which converges strongly to $(u,w).$ We use the Aubin-Lions Lemma. We can show by making use of the equation (\ref{2.6}) that $(\partial_tu_n,\partial_tw_n)$ is uniformly bounded in 
\begin{align*}
&\partial_tu^{n}\in L^{1}(0,T,B_{2,1}^{1+\frac{d}{2}-2\alpha})\cap L^{2}(0,T,B_{2,1}^{\frac{d}{2}+1-3\alpha})\;,\\&\partial_tw^{n}\in L^{1}(0,T,B_{2,1}^{1+\frac{d}{2}-2\beta})\cap L^{2}(0,T,B_{2,1}^{\frac{d}{2}+1-3\beta})\;.
\end{align*}
Since we are in this case in the whole space $\mathbb{R}^d,$ we need to combine Cantor's diagonal process with the Aubin-Lions Lemma to show that a subsequence of a weakly convergent subsequence, still denoted by $(u_n,w_n),$ has the following strongly convergent property
\begin{align*}
(u_n,w_n)\longrightarrow(u,w) \quad \text{in}\quad L^{2}(0,T,B_{2,1}^{1+\frac{d}{2}-\gamma}(Q))\;,
\end{align*}
where $\alpha\le \gamma\le 3\alpha$ and $Q\subset\mathbb{R}^d$ is a compact subset. This strong convergence property would allow us to show that $(u,w)$ is indeed a weak solution of (\ref{2.1}). This completes the proof for the existence part of Theorem \ref{thm}.
\end{proof}

\subsection{Uniqueness of weak solutions}

\begin{proof}[Proof]
Assume that $(u^{(1)},w^{(1)})$ and $(u^{(2)},w^{(2)})$ are two solutions of (\ref{2.1}) in the regularity class in (\ref{2.2}) and (\ref{2.3}). Their difference $(\widetilde{u},\widetilde{w})$ with $$\widetilde{u}=u^{(2)}-u^{(1)}\quad\text{and} \quad\widetilde{w}=w^{(2)}-w^{(1)}$$ satisfies 
\begin{align}
 \left\{
 \begin{array}{ll}
 \partial_t\widetilde{u}+(\nu+k)(-\Delta)^{\alpha}\widetilde{u}=-\mathbb{P}(u^{(2)}\cdot\nabla\widetilde{u}+\widetilde{u}\cdot\nabla u^{(1)})+2k\nabla\times\widetilde{w},\\\\
   \partial_t\widetilde{w}+\gamma(-\Delta)^{\beta}\widetilde{w}=-4k\widetilde{w}-2k\nabla\times\widetilde{u}-u^{(2)}\cdot\nabla\widetilde{w}-\widetilde{u}\cdot \nabla w^{(1)},\\\\
\nabla\cdot\widetilde{u}=0,\\\\
\widetilde{u}(x,0)=0,\;\;\widetilde{w}(x,0)=0~.
\end{array}\label{2.4}
\right.
\end{align}
 We estimate the difference $(\widetilde{u},\widetilde{w})$ in $L^2(\mathbb{R}^d)$. Dotting (\ref{2.4}) by $(\widetilde{u},\widetilde{w})$ and applying the divergence-free condition, we find
\begin{align*}
\frac{1}{2}&\frac{d}{dt}\Big(\|\widetilde{u}\|^2_{L^2}+\|\widetilde{w}\|^2_{L^2}\Big)+(\nu+k)\|\Lambda^{\alpha}\widetilde{u}\|^2_{L^2}+\gamma\|\Lambda^{\beta}\widetilde{w}\|^2_{L^2}+4k\|\widetilde{w}\|^2_{L^2}\\&=-\int u^{(2)}\cdot\nabla\widetilde{u}\cdot\widetilde{u}\,dx-\int\widetilde{u}\cdot\nabla u^{(1)}\cdot\widetilde{u}\,dx\\&\quad-\int u^{(2)}\cdot\nabla\widetilde{w}\cdot\widetilde{w}\,dx-\int\widetilde{u}\cdot\nabla w^{(1)}\cdot\widetilde{w}\,dx\\&=L_1+L_2+L_3+L_4~,
\end{align*}
where 
\begin{align*}
&L_1=-\int u^{(2)}\cdot\nabla\widetilde{u}\cdot\widetilde{u}\,dx\;,\\&L_2=-\int\widetilde{u}\cdot\nabla u^{(1)}\cdot\widetilde{u}\,dx\;,\\&L_3=-\int u^{(2)}\cdot\nabla\widetilde{w}\cdot\widetilde{w}\,dx\;,\\&L_4=-\int\widetilde{u}\cdot\nabla w^{(1)}\cdot\widetilde{w}\,dx~.
\end{align*}
Due to $\nabla\cdot u^{(2)}=0,$  we find $L_1=L_3=0$ after integration by parts. In fact,
\begin{align*}
L_1&=-\int u^{(2)}\cdot\nabla\widetilde{u}\cdot\widetilde{u}\,dx&\\&=-\int u^{(2)}\cdot\nabla(\frac{1}{2}|\widetilde{u}|^2)\,dx\\&=-\int \nabla\cdot (u^{(2)} \frac{1}{2}|\widetilde{u}|^2)\,dx\\&=0~.
\end{align*}
By H\"older's inequality and Bernstein's inequality,
\begin{align}
|L_2|&=|-\int\widetilde{u}\cdot\nabla u^{(1)}\cdot\widetilde{u}\,dx|\notag\\&\le\|\nabla u^{(1)}\|_{L^{\infty}}\|\widetilde{u}\|^2_{L^2}\notag\\&\le \sum\limits_{j\ge-1}\|\Delta_j\nabla u^{(1)}\|_{L^{\infty}}\|\widetilde{u}\|^2_{L^2}\notag\\&\le c\underbrace{\sum_{j\ge-1}2^{dj(\frac{1}{2}-\frac{1}{\infty})}2^j\|\Delta_ju^{(1)}\|_{L^2}}_{=\|u^{(1)}\|_{B_{2,1}^{1+\frac{d}{2}}}}\|\widetilde{u}\|^2_{L^2}\le c\,\|u^{(1)}\|_{B_{2,1}^{1+\frac{d}{2}}}\|\widetilde{u}\|^2_{L^2}. \label{ddff}
\end{align}
To  bound  $L_4,$ we set 
\begin{align*}
\frac{1}{p}=\frac{1}{2}-\frac{\beta}{d}\;,\quad\frac{1}{q}=\frac{\beta}{d}\;\;(or\; \frac{d}{q}=\beta)
\end{align*}
By H\"older's inequality and Bernstein's inequality,
\begin{align*}
|L_4|&=|-\int\widetilde{u}\cdot\nabla w^{(1)}\cdot\widetilde{w}\,dx|\\&\le\|\widetilde{u}\|_{L^2}\|\nabla w^{(1)}\|_{L^{q}} \|\widetilde{w}\|_{L^p}\\&\le \sum_{j\ge-1}\|\Delta_j\nabla w^{(1)}\|_{L^q}\|\widetilde{u}\|_{L^2}\|\widetilde{w}\|_{L^p}\\&\le c\sum_{j\ge-1}2^j2^{dj(\frac{1}{2}-\frac{1}{q})}\|\Delta_j w^{(1)}\|_{L^2}\|\widetilde{u}\|_{L^2}\|\widetilde{w}\|_{L^p}\\&{\le}\sum_{j\ge-1}2^{j+\frac{dj}{2}-2\beta j}\|\Delta_j w^{(1)}\|_{L^2}\|\widetilde{w}\|_{L^2}\|\widetilde{u}\|_{L^p}\\&\le c\,\|w^{(1)}\|_{B_{2,1}^{1+\frac{d}{2}-\beta}}\|\widetilde{u}\|_{L^2}
\|\Lambda^{\beta}\widetilde{w}\|_{L^2}\\&\le \frac{\gamma}{2}\|\Lambda^{\beta}\widetilde{w}\|_{L^2}^2 + c\,\|w^{(1)}\|^2_{B_{2,1}^{1+\frac{d}{2}-\beta}}\|\widetilde{u}\|^2_{L^2},
\end{align*}
where in the last inequality we have made use of
\begin{align*}
\|\widetilde{u}\|_{L^p}\le c\;\|\Lambda^{\alpha}\widetilde{u}\|_{L^2}~.
\end{align*}
Combining the estimates leads to
\begin{align}
\frac{d}{dt}\Big(\|\widetilde{u}\|^2_{L^2}&+\|\widetilde{w}\|^2_{L^2}\Big)+ 2 (\nu+k)\|\Lambda^{\alpha}\widetilde{u}\|^2_{L^2}+\gamma\,\|\Lambda^{\beta}\widetilde{w}\|^2_{L^2}+8k\,\|\widetilde{w}\|^2_{L^2}\nonumber\\&\le\Big( c\,\|u^{(1)}\|_{B_{2,1}^{1+\frac{d}{2}}}+c\,\|w^{(1)}\|^2_{B_{2,1}^{1+\frac{d}{2}-\beta}}\Big)\Big(\|\widetilde{u}\|^2_{L^2}+\|\widetilde{w}\|^2_{L^2}\Big)~.\label{4.2}
\end{align}
Since $u^{(1)}\in L^1(0,T,B_{2,1}^{1+\frac{d}{2}})$ and $w^{(1)}\in L^1(0,T,B_{2,1}^{1+\frac{d}{2}})\cap L^{\infty}(0,T, B_{2,1}^{1+\frac{d}{2}-2\beta})\;,$
\begin{align*}
&\int_0^T\|u^{(1)}(t)\|_{B_{2,1}^{1+\frac{d}{2}}}dt<\infty, \\ &\int_0^T \|w^{(1)}(t)\|^2_{B_{2,1}^{1+\frac{d}{2}-2\beta}}dt\le \int_0^T \|w^{(1)}(t)\|_{B_{2,1}^{1+\frac{d}{2}}}\, \|w^{(1)}(t)\|_{B_{2,1}^{1+\frac{d}{2}-2\beta}}dt \\
&\qquad\qquad \qquad\qquad \qquad \le  \|w^{(1)}(t)\|_{L^{\infty}(0,T,B_{2,1}^{1+\frac{d}{2}-2\beta})}\, \int_0^T \|w^{(1)}(t)\|_{B_{2,1}^{1+\frac{d}{2}}}\,dt <\infty.
\end{align*}
Applying Gronwall's inequality to (\ref{4.2}) yields
\begin{align*}
\|\widetilde{u}\|_{L^2}=\|\widetilde{w}\|_{L^2}=0\;,
\end{align*}
which leads to the desired uniqueness. This completes the proof of the uniqueness part of Theorem \ref{thm}.
\end{proof}

\vskip .3in 
\section{Proof of Theorem \ref{thm2}}

The proof of Theorem \ref{thm2} is similar to the one of Theorem \ref{thm}. To avoid repetitions, we will refer next to some inequalities already showed in the proof of Theorem \ref{thm}.
In this proof, we consider the system of equations (\ref{2.1}) with $\beta = 0$, that is
\begin{align}
 \left\{
    \begin{array}{ll}
       \partial_tu+(\nu+k)(-\Delta)^{\alpha}u+u\cdot \nabla u+\nabla\Pi-2k\nabla\times w=0,\\\\
         \partial_tw+(4k+\gamma)w +2k\nabla\times u+u\cdot \nabla w =0,\\\\   \nabla\cdot u=0,\\\\u(x,0)=u_0(x),\;\;w(x,0)=w_0(x).
 \end{array}\label{2.1a}
\right.
\end{align}
  
\subsection{Existence of a weak solution}
This subsection proves the existence part of Theorem \ref{thm2}. 
The approach is to construct a successive approximation sequence and show that the limit of a subsequence actually solves (\ref{2.1a}) in the weak sense.
\begin{proof}[Proof for the existence part of Theorem \ref{thm2}]
We consider a successive approximation $\left\{(u^{(n)},w^{(n)})\right\}$ satisfying 
\begin{align}
 \left\{
  \begin{array}{ll}
 u^{(1)}=S_2u_0\;,\quad \quad w^{(1)}=S_2w_0\\\\
  \partial_tu^{(n+1)}+(\nu+k)(-\Delta^{\alpha})u^{(n+1)}=\mathbb{P}(-u^{(n)}\cdot\nabla u^{(n+1)})+2k\nabla \times w^{(n)}\\\\
 \partial_tw^{(n+1)}=-(4k +\gamma)w^{(n+1)}-2k\nabla\times u^{(n)}-u^{(n)}\cdot\nabla w^{(n+1)}\\\\
 u^{(n+1)}(x,0)=S_{n+1}u_0\;,\quad\quad w^{(n+1)}(x,0)=S_{n+1}w_0~,
 \end{array}\label{2.66}
\right.
\end{align}
where $\mathbb{P}=I-\nabla(-\Delta)^{-1}\mathrm{div}$ is the standard Leray Projection. For \begin{align*}M=2(\|u_0\|_{B_{2,1}^{1+\frac{d}{2}-2\alpha}}+\|w_0\|_{B_{2,1}^{\frac{d}{2}}})\;,\end{align*}$T>0$ being sufficiently small and $0<\delta<1$  (to be specified later), we set 
\begin{align}
Y\equiv&\left\{(u,w)\;|\quad  \|u\|_{L^{\infty}(0,T,B_{2,1}^{1+\frac{d}{2}-2\alpha})}\le M\;,\;\;\|w\|_{L^{\infty}(0,T,B_{2,1}^{\frac{d}{2}})}\le M\;,\right.\nonumber \\
&\hphantom{M(u,v)\;|}\left.\quad\|u\|_{L^{1}(0,T,B_{2,1}^{1+\frac{d}{2}}) }\le \delta\;,\;\;\|w\|_{L^{1}(0,T,B_{2,1}^{\frac{d}{2}}) }\le \delta\right\}\;.\label{2.7}
\end{align}
We show that  $\left\{(u^{(n)},w^{(n)})\right\}$ has a subsequence that converges to the weak solution of (\ref{2.1a}). This process consists of three main steps. The first step is to show that $\left\{(u^{(n)},w^{(n)})\right\}$ is uniformly bounded in $Y$. The second step is to extract a strongly convergent subsequence via the Aulin-Lions Lemma. While the last step is to show that the limit is indeed a weak solution of (\ref{2.1a}).

\vskip .1in 
To show the uniform bound for $\left\{(u^{(n)},w^{(n)})\right\}$ in Y, we prove by induction. Clearly,
\begin{align*}
&\|u^{(1)}\|_{L^{\infty}(0,T,B_{2,1}^{1+\frac{d}{2}-2\alpha} )}=\|S_2u_0\|_{L^{\infty}(0,T,B_{2,1}^{1+\frac{d}{2}-2\alpha} )}\le M\;,\\&\|w^{(1)}\|_{L^{\infty}(0,T,B_{2,1}^{\frac{d}{2}} )}=\|S_2w_0\|_{L^{\infty}(0,T,B_{2,1}^{\frac{d}{2}} )}\le M\;.
\end{align*}
If $T>0$ is sufficiently small, then 
\begin{align*}
\|u^{(1)}\|_{L^1(0,T,B_{2,1}^{1+\frac{d}{2}} )}\le T \|S_2u_0\|_{B_{2,1}^{1+\frac{d}{2}}}\le T\,c \|u_0\|_{B_{2,1}^{1+\frac{d}{2}-2\alpha}}\le \delta\;,
\end{align*}
\begin{align*}
\|w^{(1)}\|_{L^1(0,T,B_{2,1}^{\frac{d}{2}} )}\le T \|S_2w_0\|_{B_{2,1}^{\frac{d}{2}}}\le T\,c \|w_0\|_{B_{2,1}^{\frac{d}{2}}}\le \delta~.
\end{align*}
Assuming that $(u^{(n)},w^{(n)})$ obeys the bounds defined in $Y,$ namely 
\begin{align*}
&\|u^{(n)}\|_{L^{\infty}(0,T, B_{2,1}^{1+\frac{d}{2}-2\alpha} )}\le M\;,\;\;\|w^{(n)}\|_{L^{\infty}(0,T,B_{2,1}^{\frac{d}{2}} )}\le M\;,\\&\|u^{(n)}\|_{L^{1}(0,T,B_{2,1}^{1+\frac{d}{2}}) }\le \delta\;,\;\;\|w^{(n)}\|_{L^{1}(0,T,B_{2,1}^{\frac{d}{2}}) }\le \delta\;.
\end{align*}
we prove that $\left\{(u^{(n+1)},w^{(n+1)})\right\}$ obeys the same bound for suitably selected $T>0,\; M>0$ and $\delta>0$. For sake of clarity, the proof of the four bounds is achieved in the following four steps.
\subsubsection{The estimate of $u^{(n+1)}$ in $L^\infty(0, T, B_{2,1}^{1+\frac{d}{2}-2\alpha}(\mathbb{R}^d))$}. Following the same method as in the proof of the first step of Theorem \ref{thm}, we write the inequality  
\begin{align}
\|u^{(n+1)}(t)\|_{B_{2,1}^{1+\frac{d}{2}-2\alpha}}\le \|u_0^{(n+1)}&\|_{B_{2,1}^{1+\frac{d}{2}-2\alpha}}\notag \\
& +\sum_{j\ge -1} 2^{(1+\frac{d}{2}-2\alpha)j}\int_0^{t}e^{-c_0 2^{2\alpha j}(t-\tau)}(J_1+\cdots+ J_4)\, d\tau\;,
\label{525}\end{align}
where
\begin{align*}
&J_1=c\,\|\Delta_j u^{(n+1)}\|_{L^2}\sum_{m\le j-1}2^{(1+\frac{d}{2})m}\|\Delta_m u^{(n)}\|_{L^2}\;,\\&J_2=c\,\|\Delta_j u^{(n)}\|_{L^2}\sum_{m\le j}2^{(1+\frac{d}{2})m}\|\Delta_m u^{(n+1)}\|_{L^2}\;,\\&J_3=c\,2^j\sum_{k\ge j-1} 2^{\frac{d}{2}k}\|\widetilde{\Delta}_k u^{(n+1)}\|_{L^2}\|\Delta_k u^{(n)}\|_{L^2}\;,\\& J_4=c\, 2^j \|\Delta_j w^{(n)}\|_{L^2}\;.
\end{align*}
The terms on the right hand side can be estimated as follows. Recalling the definition of $J_1$ above and using the inductive assumption on $u^{(n)},$ we have for any $t\le T\;,$
\begin{align*}
\sum_{j\ge -1} 2^{(1+\frac{d}{2}-2\alpha)j}&\int_0^t e^{c_02^{2\alpha j}(t-\tau)}J_1\,d\tau\\&=c \int_0^t\sum_{j\ge -1} 2^{(1+\frac{d}{2}-2\alpha)j}\|\Delta_j u^{(n+1)}\|_{L^2}\sum_{m\le j-1}2^{(1+\frac{d}{2})m}\|\Delta_m u^{(n)}(\tau)\|_{L^2}\,d\tau\\&\le c \int_0^t\underbrace{\sum_{j\ge -1} 2^{(1+\frac{d}{2}-2\alpha)j}\|\Delta_j u^{(n+1)}\|_{L^2}}_{=\|u^{(n+1)}\|_{B^{1+\frac{d}{2}-2\alpha}_{2,1}}}\underbrace{\sum_{j\ge -1}2^{(1+\frac{d}{2})j}\|\Delta_j u^{(n)}(\tau)\|_{L^2}}_{=\|u^{(n)}\|_{L^{1}(0,T,B^{1+\frac{d}{2}}_{2,1})}}\,d\tau\\&\le c\,\|u^{(n+1)}\|_{L^{\infty}(0,T,B^{1+\frac{d}{2}-2\alpha}_{2,1})}\|u^{(n)}\|_{L^{1}(0,T,B^{1+\frac{d}{2}}_{2,1})}\\&\le c\,\delta\,\|u^{(n+1)}\|_{L^{\infty}(0,T,B^{1+\frac{d}{2}-2\alpha}_{2,1})}~.
\end{align*}
The terms with $J_2$ and $J_3$ can be similarly estimated and obey the same bound. In fact, by Young's inequality for series convolution,
\begin{align*}
\sum_{j\le -1}&2^{(1+\frac{d}{2}-2\alpha)j}\int_0^te^{-c_02^{2\alpha j}(t-\tau)}J_2\,d\tau\\&\le c\int_0^t \sum_{j\ge -1}2^{(1+\frac{d}{2})j}\|\Delta_j u^{(n)}\|_{L^2}\sum_{m\le j}2^{2\alpha(m-j)}2^{(1+\frac{d}{2}-2\alpha)m}\|\Delta_mu^{(n+1)}(\tau)\|_{L^2}\,d\tau\\&\le c\int_0^t \|u^{(n)}(\tau)\|_{B_{2,1}^{1+\frac{d}{2}}}\|u^{(n+1)}(\tau)\|_{B_{2,1}^{1+\frac{d}{2}-2\alpha}}\,d\tau \\&\le c\,\|u^{(n)}\|_{L^1(0,T,B_{2,1}^{1+\frac{d}{2}})}\|u^{(n+1)}\|_{L^{\infty}(0,T, B_{2,1}^{1+\frac{d}{2}-2\alpha})}\\&\le c\,\delta\,\|u^{(n+1)}\|_{L^{\infty}(0,T,B_{2,1}^{1+\frac{d}{2}-2\alpha})}~.\end{align*}
Similarly the term with $J_3$ is bounded by
\begin{align*}
\sum_{j\ge -1} &2^{(1+\frac{d}{2}-2\alpha)j}\int_0^t e^{-c_02^{2\alpha j}(t-\tau)}J_3\,d\tau\\&=\int_0^t \sum_{j\ge -1}2^{(1+\frac{d}{2}-2\alpha)j}\sum_{k\ge j-1}2^{j}2^{\frac{d}{2}k}\|\widetilde{\Delta}_ku^{(n+1)}\|_{L^2}\|\Delta_k u^{(n)}\|_{L^2}\,d\tau\\&=c\int_0^t \sum_{j\ge-1}\sum_{k\ge j-1}2^{(2+\frac{d}{2}-2\alpha)(j-k)}2^{(1+\frac{d}{2})k}\|\Delta_k u^{(n)}\|_{L^2}2^{(1+\frac{d}{2}-2\alpha)k}\|\widetilde{\Delta}_ku^{(n+1)}\|_{L^2}\,d\tau\\&\le c\int_0^t \|u^{(n)}(\tau)\|_{B_{2,1}^{1+\frac{d}{2}}}\|u^{(n+1)}(\tau)\|_{B_{2,1}^{1+\frac{d}{2}-2\alpha}}\,d\tau\\&\le c\,\|u^{(n)}\|_{L^1(0,T, B_{2,1}^{1+\frac{d}{2}})}\|u^{(n+1)}\|_{L^{\infty}(0,T, B_{2,1}^{1+\frac{d}{2}-2\alpha})}\\&\le c\,\delta\, \|u^{(n+1)}\|_{L^{\infty}(0,T, B_{2,1}^{1+\frac{d}{2}-2\alpha})}~.\end{align*}
Now for the term with $J_4$ we write
\begin{align*}
\sum_{j\ge -1}2^{(1+\frac{d}{2}-2\alpha)j}\int_0^t e^{-c_02^{2\alpha j(t-\tau)}}J_4\,d\tau&=\sum_{j\ge -1}2^{(1+\frac{d}{2}-2\alpha)j}\int_0^t e^{-c_02^{2\alpha j(t-\tau)}} c 2^j\|\Delta_j w^{(n)}\|_{L^2}\,d\tau\\&\le \sum_{j\ge -1}\int_0^t  c\,2^{(2+\frac{d}{2}-2\alpha)j}\|\Delta_j w^{(n)}\|_{L^2}\,d\tau\\&\underbrace{\le}_{\text{since}\;\alpha\ge 1}c\,\|w^{(n)}\|_{L^1(0,T,B_{2,1}^{\frac{d}{2}})}\\&\le c\,\delta~.
\end{align*}
Collecting the bounds above and inserting them in (\ref{525}), we find for any $t\le T$ 
\begin{align*}
\|u^{(n+1)}(t)\|_{B_{2,1}^{1+\frac{d}{2}-2\alpha}}\le \|u_0^{(n+1)}\|_{B_{2,1}^{1+\frac{d}{2}-2\alpha}}+c\,\delta\,\|u^{(n+1)}\|_{L^{\infty}(0,T,B^{1+\frac{d}{2}-2\alpha}_{2,1})}+c\,\delta~.
\end{align*}
Therefore
\begin{align*}
\|u^{(n+1)}(t)\|_{L^{\infty}(0,T,B_{2,1}^{1+\frac{d}{2}-2\alpha})}\le \|u_0^{(n+1)}\|_{B_{2,1}^{1+\frac{d}{2}-2\alpha}}+c\,\delta\|u^{(n+1)}\|_{L^{\infty}(0,T,B_{2,1}^{1+\frac{d}{2}-2\alpha})}+c\,\delta~.
\end{align*}
Choosing $\delta$ such that $c\,\delta \le \min(\frac{1}{4},\frac{M}{4})$ we get
\begin{align*}
\|u^{(n+1)}(t)\|_{L^{\infty}(0,T,B_{2,1}^{1+\frac{d}{2}-2\alpha})}\le \frac{M}{2}+\frac{1}{4}\|u^{(n+1)}\|_{L^{\infty}(0,T,B_{2,1}^{1+\frac{d}{2}-2\alpha})}+\frac{M}{4}~,
\end{align*}
which implies
\begin{align*}
\|u^{(n+1)}(t)\|_{L^{\infty}(0,T,B_{2,1}^{1+\frac{d}{2}-2\alpha})}\le M~.
\end{align*}
\subsubsection{The estimate of $\|u^{(n+1)}\|_{L^1(0,T,B^{1+\frac{d}{2}})}$}. Following the same method as in the proof of the second step of Theorem \ref{thm}, we write the inequality 
\begin{align}
\|u^{(n+1)}\|_{L^1(0,T,B_{2,1}^{1+\frac{d}{2}})}\le \int_0^T\sum_{j\ge -1}&2^{(1+\frac{d}{2})j}e^{-c_02^{2\alpha j}t}\|\Delta_j u_0^{(n+1)}\|_{L^2}\,dt\nonumber\\&+\int_0^T\int_0^s\sum_{j\ge -1}2^{(1+\frac{d}{2})j}e^{-c_02^{2\alpha j}(s-\tau)}(J_1+\cdots +J_4)\,d\tau \,ds\;.
\label{524}\end{align}
We estimate the terms on the right and start with the first term.
\begin{align*}
\int_0^T\sum_{j\ge -1}2^{(1+\frac{d}{2})j}e^{-c_02^{2\alpha j}t}\|\Delta_ju_0^{(n+1)}\|_{L^2}\,dt&=c\sum_{j\ge -1}2^{(1+\frac{d}{2}-2\alpha)j}(1-e^{-c_02^{2\alpha j}T})\|\Delta_j u_0^{(n+1)}\|_{L^2}\;.
\end{align*}
Since $u_0\in B_{2,1}^{1+\frac{d}{2}-2\alpha}\;,$ then by the Dominated Convergence Theorem
\begin{align*}
\lim_{T\to 0}\sum_{j\ge -1}2^{(1+\frac{d}{2}-2\alpha)j}(1-e^{-c_02^{2\alpha j}T})\|\Delta_j u_0^{(n+1)}\|_{L^2}=0\;.
\end{align*}
Therefore, we can choose $T$ sufficiently small such that
\begin{align*}
\int_0^T\sum_{j\ge -1}2^{(1+\frac{d}{2})j}e^{-c_02^{2\alpha j}t}\|\Delta_j u_0^{(n+1)}\|_{L^2}\le \frac{\delta}{4}\;.
\end{align*}
Applying Young's inequality for the time convolution, we have 
\begin{align*}
\int_0^T&\sum_{j\ge -1}2^{(1+\frac{d}{2})j}\int_0^se^{-c_02^{2\alpha j}(s-\tau)} J_1\,d\tau\,ds\\&=c\int_0^T\sum_{j\ge -1}2^{(1+\frac{d}{2})j}\int_0^se^{-c_02^{2\alpha j}(s-\tau)} \|\Delta_j u^{(n+1)}(\tau)\|_{L^2}\sum_{m\le j -1}2^{(1+\frac{d}{2})m} \|\Delta_m u^{(n)}\|_{L^2}\,d\tau\,ds\\&\le c\sum_{j\ge -1}2^{(1+\frac{d}{2})j}\int_0^T\|\Delta_j u^{(n+1)}(\tau)\|_{L^2}\sum_{m\le j -1}2^{(1+\frac{d}{2})m} \|\Delta_m u^{(n)}\|_{L^2}\,d\tau \int_0^Te^{-c_02^{2\alpha j}s}\,ds~.
\end{align*}
Then, using the fact that there exists $c_2> 0$ satisfying for $j\ge 0\;,$
\begin{align}
\int_0^Te^{-c_02^{2\alpha j}s}ds\le c\,2^{-2\alpha j}(1-e^{-c_2T})\;,\label{sou2}
\end{align}
we get
\begin{align*}
\int_0^T\sum_{j\ge -1}2^{(1+\frac{d}{2})j}&\int_0^se^{-c_02^{2\alpha j}(s-\tau)} J_1\,d\tau\,ds\\&\le c\,(1-e^{-c_2T})\int_0^T\sum_{j\ge -1}2^{(1+\frac{d}{2}-2\alpha)j}\|\Delta_j u^{(n+1)}(\tau)\|_{L^2}\sum_{m\le j -1}2^{(1+\frac{d}{2})m} \|\Delta_m u^{(m)}\|_{L^2}\,d\tau\\&\le c\,(1-e^{-c_2T})\int_0^T\sum_{j\ge -1}2^{(1+\frac{d}{2}-2\alpha)j}\|\Delta_j u^{(n+1)}(\tau)\|_{L^2}\sum_{j\ge-1}2^{(1+\frac{d}{2})j} \|\Delta_j u^{(n)}\|_{L^2}\,d\tau\\&\le c\,(1-e^{-c_2T})\|u^{(n+1)}\|_{L^{\infty}(0,T,B_{2,1}^{1+\frac{d}{2}-2\alpha})}\underbrace{\|u^{(n)}\|_{L^{1}(0,T,B_{2,1}^{1+\frac{d}{2}})}}_{\le \delta}\\&\le c\,\delta\,(1-e^{-c_2T})\|u^{(n+1)}\|_{L^{\infty}(0,T,B_{2,1}^{1+\frac{d}{2}-2\alpha})}\;, 
\end{align*}
The terms with $J_2$ and $J_3$ can be similarly estimated and obey the same bound.
\begin{align*}\int_0^T&\sum\limits_{j\ge -1}2^{(1+\frac{d}{2})j}\int_0^se^{-c_02^{2\alpha j}(s-\tau)} J_2\,d\tau\,ds\\&=c\int_0^T\sum\limits_{j\ge -1}2^{(1+\frac{d}{2})j}\int_0^se^{-c_02^{2\alpha j}(s-\tau)} \|\Delta_ju^{(n)}\|_{L^2}\sum\limits_{m\le j}2^{(1+\frac{d}{2})m}\|\Delta_mu^{(n+1)}\|_{L^2}d\tau\,ds\\&\le c\sum\limits_{j\ge -1}2^{(1+\frac{d}{2})j}\int_0^T\|\Delta_ju^{(n)}\|_{L^2}\sum\limits_{m\le j}2^{(1+\frac{d}{2})m}\|\Delta_mu^{(n+1)}\|_{L^2}d\tau\int_0^Te^{-c_02^{2\alpha j}s}\,ds~.
\end{align*}
Hence due to (\ref{sou2})\,
\begin{align*}
\int_0^T&\sum\limits_{j\ge -1}2^{(1+\frac{d}{2})j}\int_0^se^{-c_02^{2\alpha j}(s-\tau)} J_2\,d\tau\,ds\\&\le c (1-e^{-c_2T})\int_0^T\sum\limits_{j\ge -1}2^{(1+\frac{d}{2}-2\alpha)j}\|\Delta_ju^{(n)}\|_{L^2}\sum\limits_{m\le j-1}2^{(1+\frac{d}{2})m}\|\Delta_mu^{(n+1)}\|_{L^2}d\tau\\&\le c (1-e^{-c_2T})\int_0^T\sum\limits_{j\ge -1}2^{(1+\frac{d}{2}-2\alpha)j}\|\Delta_ju^{(n+1)}\|_{L^2}\sum\limits_{j\ge -1}2^{(1+\frac{d}{2})j}\|\Delta_ju^{(n)}\|_{L^2}d\tau \\&\le c(1-e^{-c_2T})\|u^{(n+1)}\|_{L^{\infty}(0,T,B_{2,1}^{1+\frac{d}{2}-2\alpha})}\|u^{(n)}\|_{L^{1}(0,T,B_{2,1}^{1+\frac{d}{2}})}\\&\le c(1-e^{-c_2T})\,\delta\,\|u^{(n+1)}\|_{L^{\infty}(0,T,B_{2,1}^{1+\frac{d}{2}-2\alpha})}\;.\end{align*}
Similarly
\begin{align*}\int_0^T\sum\limits_{j\ge -1}2^{(1+\frac{d}{2})j}&\int_0^se^{-c_02^{2\alpha j}(s-\tau)} J_3\,d\tau\,ds
\\&=c\int_0^T\sum\limits_{j\ge -1}2^{(2+\frac{d}{2})j}\int_0^se^{-c_02^{2\alpha j}(s-\tau)}\sum\limits_{k\ge j-1}2^{\frac{d}{2}k}\|\Delta_ku^{(n)}\|_{L^2}\|\widetilde{\Delta}_ku^{(n+1)}\|_{L^2}\,d\tau\,ds
\\&\le c\sum\limits_{j\ge -1}2^{(2+\frac{d}{2})j}\int_0^T\sum\limits_{k\ge j-1}2^{\frac{d}{2}k}\|\Delta_ku^{(n)}\|_{L^2} \|\widetilde{\Delta}_ku^{(n+1)}\|_{L^2}d\tau\int_0^Te^{-c_02^{2\alpha j}s}\,ds~.
\end{align*}
By (\ref{sou2}) and the inequality above, 
\begin{align*}
\int_0^T\sum\limits_{j\ge -1}2^{(1+\frac{d}{2})j}&\int_0^se^{-c_02^{2\alpha j}(s-\tau)} J_3\,d\tau\,ds
\\&\le c\,(1-e^{-c_2T})\int_0^T\sum\limits_{j\ge -1}2^{(2+\frac{d}{2}-2\alpha)j}\sum\limits_{k\ge j-1}2^{\frac{d}{2}k}\|\Delta_ku^{(n)}\|_{L^2}\|\widetilde{\Delta}_ku^{(n+1)}\|_{L^2}\,d\tau
\\&\le c\,(1-e^{-c_2T})\int_0^T\sum\limits_{j\ge -1}2^{(2+\frac{d}{2}-2\alpha)j}\|\widetilde{\Delta}_ju^{(n+1)}\|_{L^2}\sum\limits_{j\ge -1}2^{(\frac{d}{2})j}\|\Delta_ju^{(n)}\|_{L^2}\,d\tau
\\&\le c\,(1-e^{-c_2T})\int_0^T\sum\limits_{j\ge -1}2^{(1+\frac{d}{2}-2\alpha)j}\|\widetilde{\Delta}_ju^{(n+1)}\|_{L^2}\sum\limits_{j\ge -1}2^{(1+\frac{d}{2})j}\|\Delta_ju^{(n)}\|_{L^2}\,d\tau 
\\&\le c\,(1-e^{-c_2T})\|u^{(n+1)}\|_{L^{\infty}(0,T,B_{2,1}^{1+\frac{d}{2}-2\alpha})}\|u^{(n)}\|_{L^{1}(0,T,B_{2,1}^{1+\frac{d}{2}})}\\&\le c\,(1-e^{-c_2T})\,\delta\,\|u^{(n+1)}\|_{L^{\infty}(0,T,B_{2,1}^{1+\frac{d}{2}-2\alpha})}\;.\end{align*}
The term with $J_4$ is bounded by\begin{align*}\int_0^T\sum\limits_{j\ge -1}2^{(1+\frac{d}{2})j}&\int_0^se^{-c_02^{2\alpha j}(s-\tau)} J_4\,d\tau\,ds
\\&=c\int_0^T\sum\limits_{j\ge -1}2^{(2+\frac{d}{2})j}\int_0^se^{-c_02^{2\alpha j}(s-\tau)} \|\Delta_jw^{(n)}\|_{L^2}\,d\tau\,ds
\\&\le c\sum\limits_{j\ge -1}2^{(2+\frac{d}{2})j}\int_0^T\|\Delta_jw^{(n)}\|_{L^2}d\tau\int_0^Te^{-c_02^{2\alpha j}s}\,ds
\\&\le c\,(1-e^{-c_2T})\int_0^T\sum\limits_{j\ge -1}2^{(2+\frac{d}{2}-2\alpha)}\|\Delta_jw^{(n)}\|_{L^2}\,d\tau
\\&\underset{\text{since}\,\alpha\ge 1}{\le} c\,(1-e^{-c_2T})\int_0^T\sum\limits_{j\ge -1}2^{(\frac{d}{2})j}\|\Delta_jw^{(n)}\|_{L^2}\,d\tau
\\&=c (1-e^{-c_2T})\|w^{(n)}\|_{L^{1}(0,T,B_{2,1}^{\frac{d}{2}})}\;.\end{align*}
Collecting the estimates above and inserting them in (\ref{524}) leads to
\begin{align*}
\|u^{(n+1)}\|_{L^{1}(0,T,B_{2,1}^{1+\frac{d}{2}})}&\le \frac{\delta}{4}+c\,\delta(1-e^{-c_2T})\|u^{(n+1)}\|_{L^{\infty}(0,T,B_{2,1}^{1+\frac{d}{2}-2\alpha})}\\&\quad\quad+c\,(1-e^{-c_2T})\|w^{(n)}\|_{L^{1}(0,T,B_{2,1}^{\frac{d}{2}})}\\&\le \frac{\delta}{4}+c\,\delta(1-e^{-c_2T})M+c\,(1-e^{-c_2T})\delta\;.
\end{align*}
Choosing T sufficiently small such that 
$c(1-e^{-c_2T})\le \min(\frac{1}{4M} , \frac{1}{2})$ we get
\begin{align*}
\|u^{(n+1)}\|_{L^{1}(0,T,B_{2,1}^{1+\frac{d}{2}})}\le \frac{\delta}{4}+\frac{\delta}{4}+\frac{\delta}{2}=\delta\;.
\end{align*}

\subsubsection{The estimate of $w^{(n+1)}$ in $L^\infty(0, T, B_{2,1}^{\frac{d}{2}}(\mathbb{R}^d)).$}  Applying $\Delta_j$ to the third equation in (\ref{2.66})  and then dotting with $\Delta_jw^{(n+1)},$ we obtain
\begin{align}
\frac{1}{2}\frac{d}{dt}\|\Delta_jw^{(n+1)}\|_{L^2}^2 +(4k+\gamma)\|\Delta_jw^{(n+1}\|^2_{L^2}&=-2k\int \Delta_j(\nabla\times u^{(n)})\Delta_jw^{(n+1)}dx\nonumber\\&\quad-\int \Delta_j( u^{(n)}\cdot\nabla w^{(n+1)}) \Delta_jw^{(n+1)}dx\nonumber\\&=B_1+B_2\;,
\label{3.10}\end{align}
where 
\begin{align*}
&B_1=-2k\int \Delta_j(\nabla\times u^{(n)})\Delta_jw^{(n+1)}dx\;,\\
&B_2=-\int \Delta_j( u^{(n)}\cdot\nabla w^{(n+1)})\Delta_jw^{(n+1)}dx\;.
\end{align*}
By H\"older's inequality and Bernstein's inequality
\begin{align*}
|B_1|&=|-2k\int \Delta_j(\nabla\times u^{(n)})\Delta_jw^{(n+1)}dx|\\&\le 2k\,\|\Delta_j(\nabla\times u^{(n)})\|_{L^2}\|\Delta_jw^{(n+1)}dx\|_{L^2}\\&\le c\,2^j\,\|\Delta_j u^{(n)}\|_{L^2}\|\Delta_jw^{(n+1)}\|_{L^2}\;.
\end{align*}
By Lemma \ref{lem2.3},
\begin{align*}
|B_2|&=|-\int \Delta_j( u^{(n)}\cdot\nabla w^{(n+1)})\Delta_jw^{(n+1)}\,dx|\\&\le c\,\|\Delta_j w^{(n+1)}\|_{L^2}^2\sum\limits_{m\le j-1}2^{(1+\frac{d}{2})m}\|\Delta_mu^{(n)}\|_{L^2}\\&+c\,\|\Delta_j w^{(n+1)}\|_{L^2}\|\Delta_j u^{(n)}\|_{L^2}\sum\limits_{m\le j}2^{(1+\frac{d}{2})m}\|\Delta_mw^{(n+1)}\|_{L^2}\\&+c\,\|\Delta_j w^{(n+1)}\|_{L^2}2^j\sum\limits_{k\le j-1}2^{\frac{d}{2}k}\|\widetilde{\Delta}_k w^{(n+1)}\|_{L^2}\|\Delta_ku^{(n)}\|_{L^2}\;.
\end{align*}
Inserting the estimates above in (\ref{3.10}) and eliminating $\|\Delta_j w^{(n+1)}\|_{L^2}$ from both sides of the inequality, we obtain
\begin{align}
\frac{d}{dt}\|\Delta_jw^{(n+1)}\|_{L^2}+(8k+2\gamma)\|\Delta_jw^{(n+1)}\|_{L^2}\le & c\,2^j \|\Delta_ju^{(n)}\|_{L^2}\notag\\
& +c\,\|\Delta_jw^{(n+1)}\|_{L^2}\sum\limits_{m\le j-1}2^{(1+\frac{d}{2})m}\|\Delta_mu^{(n)}\|_{L^2}\nonumber\\&+c\,\|\Delta_j u^{(n)}\|_{L^2}\sum\limits_{m\le j}2^{(1+\frac{d}{2})m}\|\Delta_mw^{(n+1)}\|_{L^2}\nonumber\\&+c\,2^j\sum\limits_{k\ge j-1}2^{\frac{d}{2}k}\|\widetilde{\Delta}_k w^{(n+1)}\|_{L^2}\|\Delta_ku^{(n)}\|_{L^2}~.
\label{3.11}\end{align}
Integrating (\ref{3.11}) in time yields, for any $t\le T$,
\begin{eqnarray}
\|\Delta_jw^{(n+1)}\|_{L^2}\le e^{-(8k+2\gamma)t}\|\Delta_jw_0^{(n+1)}\|_{L^2}+\int_0^te^{-(8k+2\gamma)(t-\tau)}(K_1+\cdots+K_4)d\tau\label{3.125}~,
\end{eqnarray}
where 
\begin{align*}
&K_1=c\,2^j \|\Delta_ju^{(n)}\|\;,\\
&K_2=c\,\|\Delta_jw^{(n+1)}\|_{L^2}\sum\limits_{m\le j-1}2^{(1+\frac{d}{2})m}\|\Delta_mu^{(n)}\|_{L^2}\;,\\
&K_3=c\,\|\Delta_j u^{(n)}\|_{L^2}\sum\limits_{m\le j}2^{(1+\frac{d}{2})m}\|\Delta_mw^{(n+1)}\|_{L^2}\;,\\
&K_4=c\,2^j\sum\limits_{k\ge j-1}2^{\frac{d}{2}k}\|\widetilde{\Delta}_k w^{(n+1)}\|_{L^2}\|\Delta_ku^{(n)}\|_{L^2}\;.
\end{align*}
We multiply (\ref{3.125}) by $2^{(\frac{d}{2}) j}$ and sum over $j$ to get 
\begin{align}
\|w^{(n+1)}\|_{B_{2,1}^{\frac{d}{2}}}\le \|w^{(n+1)}_0\|_{B_{2,1}^{\frac{d}{2}}}+\sum_{j\ge -1}\int_0^te^{-(8k+2\gamma)(t-\tau)}2^{(\frac{d}{2})j}(K_1+\cdots+K_4)\,d\tau\;.
\label{3.122}\end{align}
The term with $K_1$ is bounded by
\begin{align*}
\sum_{j\ge -1}2^{\frac{d}{2}j}\int_0^te^{-(8k+2\gamma)(t-\tau)}K_1\,d\tau&\le \int_0^t\sum_{j\ge -1} c\,2^{(1+\frac{d}{2})j} \|\Delta_ju^{(n)}(\tau)\|_{L^2}\,d\tau\\&\le c\,\|u^{(n)}\|_{L^{1}(0,T,B_{2,1}^{1+\frac{d}{2}})}\;.
\end{align*}
The terms with $K_2$ through $K_4$ can be bounded suitably and obey the same bound. In fact, for the term involving $K_2$ we write 
\begin{align*}
\sum_{j\ge -1}&2^{(\frac{d}{2})j}\int_0^te^{-(8k+2\gamma)(t-\tau)}c\,\|\Delta_jw^{(n+1)}(\tau)\|_{L^2}\sum_{m\le j-1}2^{(1+\frac{d}{2})m}\|\Delta_mu^{(n)}(\tau)\|_{L^2}\,d\tau
\\&\le c\int_0^t\sum_{j\ge -1}2^{(\frac{d}{2})j}\|\Delta_jw^{(n+1)}(\tau)\|_{L^2}\sum_{m\le j-1}2^{(1+\frac{d}{2})m}\|\Delta_mu^{(n)}(\tau)\|_{L^2}d\tau
\\&\le c\int_0^t\|w^{(n+1)}\|_{B_{2,1}^{\frac{d}{2}}}\|u^{(n)}\|_{B_{2,1}^{1+\frac{d}{2}}}\,d\tau
\\&\le c\,\|w^{(n+1)}\|_{L^{\infty}(0,T,B_{2,1}^{\frac{d}{2}})}\|u^{(n)}\|_{L^{1}(0,T,B_{2,1}^{1+\frac{d}{2}})}\;.
\end{align*}
Similarly the term with $K_3$ is bounded by 
\begin{align*}
\sum_{j\ge -1}&2^{(\frac{d}{2})j}\int_0^te^{-(8k+2\gamma)(t-\tau)}c\,\|\Delta_ju^{(n)}\|_{L^2} \sum_{m\le j}2^{(1+\frac{d}{2})m}\|\Delta_mw^{(n+1)}\|_{L^2}d\tau\\&\le c\int_0^t\sum_{j\ge -1}2^{(\frac{d}{2})j}\|\Delta_ju^{(n)}\|_{L^2} \sum_{j\ge -1}2^{(1+\frac{d}{2})j}\|\Delta_jw^{(n+1)}\|_{L^2}\,d\tau\\&\le c\,\|w^{(n+1)}\|_{L^{\infty}(0,T,B_{2,1}^{\frac{d}{2}})}\|u^{(n)}\|_{L^{1}(0,T,B_{2,1}^{1+\frac{d}{2}})}\;.
\end{align*}
For the term with $K_4$ we write
\begin{align*}
\sum_{j\ge -1}&2^{(\frac{d}{2})j}\int_0^tc\,2^j\sum_{k\ge j-1}\|\widetilde{\Delta}_kw^{(n+1)}\|_{L^2} \|\Delta_ku^{(n)}\|_{L^2}2^{\frac{dk}{2}}\,d\tau\\&\le c\int_0^t\sum_{j\ge -1}2^{(1+\frac{d}{2})j}\|\Delta_ju^{(n)}\|_{L^2} \sum_{j\ge -1}2^{(\frac{d}{2})j}\|\widetilde{\Delta}_jw^{(n+1)}\|_{L^2}d\tau \\&\le c\,\|w^{(n+1)}\|_{L^{\infty}(0,T,B_{2,1}^{\frac{d}{2}})}\|u^{(n)}\|_{L^{1}(0,T,B_{2,1}^{1+\frac{d}{2}})}\;.
\end{align*}
Collecting the estimates and inserting them in (\ref{3.122}), we obtain for any $t\le T$
\begin{align*}
\|w^{(n+1)}(t)\|_{B_{2,1}^{\frac{d}{2}}}&\le \|w_0^{(n+1)}\|_{B_{2,1}^{\frac{d}{2}}}+c\,\|u^{(n)}\|_{L^{1}(0,T,B_{2,1}^{1+\frac{d}{2}})}+c\,\|w^{(n+1)}\|_{L^{\infty}(0,T,B_{2,1}^{\frac{d}{2}})}\|u^{(n)}\|_{L^{1}(0,T,B_{2,1}^{1+\frac{d}{2}})}\\&\le \|w_0^{(n+1)}\|_{B_{2,1}^{\frac{d}{2}}}+c\,\delta+c \,\delta\,\|w^{(n+1)}\|_{L^{\infty}(0,T,B_{2,1}^{\frac{d}{2}})}~.
\end{align*}
Therefore 
\begin{align*}
\|w^{(n+1)}(t)\|_{L^{\infty}(0,T,B_{2,1}^{\frac{d}{2}})}\le \|w_0^{(n+1)}\|_{B_{2,1}^{\frac{d}{2}}}+c\,\delta+c \,\delta\,\|w^{(n+1)}\|_{L^{\infty}(0,T,B_{2,1}^{\frac{d}{2}})}~.
\end{align*}
Choosing $c\,\delta\le\min(\frac{1}{4},\frac{M}{4})$ we get
\begin{align*}
\|w^{(n+1)}(t)\|_{L^{\infty}(0,T,B_{2,1}^{\frac{d}{2}})}\le\frac{M}{2}+\frac{M}{4}+\frac{1}{4}\|w^{(n+1)}(t)\|_{L^{\infty}(0,T,B_{2,1}^{\frac{d}{2}})}~,
\end{align*}
which implies
\begin{align*}
\|w^{(n+1)}(t)\|_{L^{\infty}(0,T,B_{2,1}^{\frac{d}{2}})}\le M~.
\end{align*}

\subsubsection{The estimate of $\|w^{(n+1)}(t)\|_{L^{1}(0,T,B_{2,1}^{\frac{d}{2}})}.$} We multiply (\ref{3.125}) by $2^{(\frac{d}{2})j},$ sum over $j$ and integrate in time to get 
\begin{align}
\|w^{(n+1)}\|_{L^1(0,T,B_{2,1}^{\frac{d}{2}})}&\le\int_0^T\sum_{j\ge -1}2^{\frac{d}{2}j}e^{-(8k+2\gamma)t} \|\Delta_jw^{(n+1)}_0\|_{L^2}\,dt\nonumber\\&\quad\quad+ \int_0^T\sum_{j\ge -1}2^{\frac{d}{2}j}\int_0^se^{-(8k+2\gamma)(s-\tau)}(K_1+\cdots+ K_4)\,d\tau\,ds~.
\label{3.133}\end{align}
Clearly 
\begin{align*}
\int_0^T\sum_{j\ge -1}2^{\frac{d}{2}j}e^{-(8k+2\gamma)t} \|\Delta_jw^{(n+1)}_0\|_{L^2}dt=c\sum_{j\ge -1}2^{\frac{d}{2}j}(1-e^{-(8k+2\gamma)T})\|\Delta_jw^{(n+1)}_0\|_{L^2}\;.
\end{align*}
Since $w_0\in B_{2,1}^{\frac{d}{2}},$ it follows from the Dominated Convergence Theorem that
\begin{align*}
\lim_{T\to 0}\sum_{j\ge -1}2^{\frac{d}{2}j}(1-e^{-(8k+2\gamma)T})\|\Delta_jw^{(n+1)}_0\|_{L^2}=0\;.
\end{align*}
Therefore, we can choose $T$ sufficiently small such that
\begin{align*}
\int_0^T\sum_{j\ge -1}2^{\frac{d}{2}j}e^{-(8k+2\gamma)t}\|\Delta_jw^{(n+1)}_0\|_{L^2}\, dt\le \frac{\delta}{2}\;.
\end{align*}
Applying the Young's inequality for the time convolution, the term with $K_1$ is bounded by 
\begin{align*}
\int_0^T\sum_{j\ge -1}2^{\frac{d}{2}j}&\int_0^se^{-(8k+2\gamma)(s-\tau)} K_1\,d\tau \,ds\\&=c\int_0^T\sum_{j\ge -1}2^{(1+\frac{d}{2})j}\int_0^se^{-(8k+2\gamma)(s-\tau)} \|\Delta_ju^{(n)}\|_{L^2}\,d\tau\,ds
\\&\le \Big(c\sum_{j\ge -1}2^{(1+\frac{d}{2})j}\int_0^T\|\Delta_ju^{(n)}\|_{L^2}\,d\tau\Big)\Big(\int_0^Te^{-(8k+2\gamma)s}ds\Big)
\\&\le c\,(1-e^{-(8k+2\gamma)T})\underbrace{\int_0^T\sum_{j\ge -1}2^{(1+\frac{d}{2})j}\|\Delta_ju^{(n)}\|_{L^2}\,d\tau}_{=\|u^{(n)}\|_{L^{1}(0,T,B_{2,1}^{1+\frac{d}{2}})}\le \delta}\\&\le c\,(1-e^{-(8k+2\gamma)T})\delta~.
\end{align*}
Similarly by applying Young's inequality for the time convolution the term with $K_2$ is bounded by 
\begin{align*}
\int_0^T\sum_{j\ge -1}2^{\frac{d}{2}j}&\int_0^se^{-(8k+2\gamma)(s-\tau)} K_2\,d\tau \,ds\\&= c\int_0^T\sum_{j\ge -1}2^{\frac{d}{2}j}\int_0^se^{-(8k+2\gamma)(s-\tau)} \|\Delta_jw^{(n+1)}\|_{L^2}\sum_{m\le j-1}2^{(1+\frac{d}{2})m}\|\Delta_mu^{(n)}\|_{L^2} d\tau ds\\&\le \Big(c\sum_{j\le -1}2^{\frac{d}{2}j}\int_0^T\|\Delta_jw^{(n+1)}\|_{L^2} \sum_{m\le j-1}2^{(1+\frac{d}{2})m}\|\Delta_mu^{(n)}\|_{L^2}\,d\tau\Big)  \Big(\int_0^Te^{-(8k+2\gamma)s}ds\Big)\\&\le 
 c\Big(1-e^{-(8k+2\gamma)T})\int_0^T\sum_{j\le -1}2^{\frac{d}{2}j}\|\Delta_jw^{(n+1)}\|_{L^2} \sum_{m\le j-1}2^{(1+\frac{d}{2})m}\|\Delta_mu^{(n)}\|_{L^2}\,d\tau  
\\&\le c\,(1-e^{-(8k+2\gamma)T})\int_0^T\sum_{j\ge -1}2^{\frac{d}{2}j}\|\Delta_jw^{(n+1)}\|_{L^2}\sum_{j\ge -1}2^{(1+\frac{d}{2})j}\|\Delta_ju^{(n)}\|_{L^2}\,d\tau\\&\le c\,(1-e^{-(8k+2\gamma)T})\|w^{(n+1)}\|_{L^{\infty}(0,T,B_{2,1}^{\frac{d}{2}})}\|u^{(n)}\|_{L^{1}(0,T,B_{2,1}^{1+\frac{d}{2}})}\\&\le c(1-e^{-(8k+2\gamma)T})\delta M\;.
\end{align*}
Applying Young's inequality for the time convolution, the term with $K_3$ is bounded by
\begin{align*}
&\int_0^T\sum_{j\ge -1}2^{\frac{d}{2}j}\int_0^se^{-(8k+2\gamma)(s-\tau)}K_3\,d\tau ds\\&=c\Big(\int_0^T\sum_{j\ge -1}2^{\frac{d}{2}j}\int_0^se^{-(8k+2\gamma)(s-\tau)}  \|\Delta_ju^{(n)}\|_{L^2}\Big)\Big(\sum_{m\le j}2^{(1+\frac{d}{2})m}\|\Delta_mw^{(n+1)}\|_{L^2}\Big)\,d\tau\,ds\\&\le\Big( c\sum_{j\ge -1}2^{\frac{d}{2}j}\int_0^T \|\Delta_ju^{(n)}\|_{L^2}\sum_{m\le j}2^{(1+\frac{d}{2})m}\|\Delta_mw^{(n+1)}\|_{L^2}\,d\tau\Big)\Big(\int_0^Te^{-(8k+2\gamma)s} ds\Big)
\\&\le c\,(1-e^{-(8k+2\gamma)T})\int_0^T\sum_{j\ge -1}2^{\frac{d}{2}j}\|\Delta_ju^{(n)}\|_{L^2}\sum_{m\le j}2^{(1+\frac{d}{2})m}\|\Delta_mw^{(n+1)}\|_{L^2}\,d\tau
\\&\le c\,(1-e^{-(8k+2\gamma)T})\int_0^T\sum_{j\ge -1}2^{(1+\frac{d}{2})j}\|\Delta_ju^{(n)}\|_{L^2}\sum_{j\ge -1}2^{\frac{d}{2}j}\|\Delta_jw^{(n+1)}\|_{L^2}\,d\tau
\\&\le c(1-e^{-(8k+2\gamma)T})\|w^{(n+1)}\|_{L^{\infty}(0,T,B_{2,1}^{\frac{d}{2}})}\|u^{(n)}\|_{L^{1}(0,T,B_{2,1}^{1+\frac{d}{2}})}\\&\le c\,(1-e^{-(8k+2\gamma)T})\delta M\;.
\end{align*}
Similarly the term with $K_4$ is bounded by 
\begin{align*}
\int_0^T\sum_{j\ge -1}2^{(\frac{d}{2})j}&\int_0^se^{-(8k+2\gamma)(s-\tau)}K_4\,d\tau ds
\\&=c\int_0^T\sum_{j\ge -1}2^{(1+\frac{d}{2})j}\int_0^se^{-(8k+2\gamma)(s-\tau)}\sum_{k\ge j-1}2^{\frac{d}{2}k}\|\widetilde{\Delta}_kw^{(n+1)}\|_{L^2}\|\Delta_ku^{(n)}\|_{L^2}d\tau ds
\\&\le \Big(c\sum_{j\ge -1}2^{(1+\frac{d}{2})j}\int_0^T\sum_{k\ge j-1}  2^{\frac{d}{2}k}\|\widetilde{\Delta}_kw^{(n+1)}\|_{L^2}\|\Delta_ku^{(n)}\|_{L^2}d\tau\Big)\Big(\int_0^Te^{-(8k+2\gamma)s}ds\Big)
\\&\le c(1-e^{-(8k+2\gamma)T})\int_0^T\sum_{j\ge -1}2^{(1+\frac{d}{2})j}\sum_{k\ge j-1}2^{\frac{d}{2}k}\|\widetilde{\Delta}_kw^{(n+1)}\|_{L^2}\|\Delta_ku^{(n)}\|_{L^2}d\tau
\\&\le c(1-e^{-(8k+2\gamma)T})\int_0^T\sum_{j\ge -1}2^{(1+\frac{d}{2})j}\|\Delta_ju^{(n)}\|_{L^2}\sum_{j\ge -1}2^{\frac{d}{2}j}\|\widetilde{\Delta}_jw^{(n+1)}\|_{L^2}d\tau
\\&\le c(1-e^{-(8k+2\gamma)T})\|w^{(n+1)}\|_{L^{\infty}(0,T,B_{2,1}^{\frac{d}{2}})}\|u^{(n)}\|_{L^{1}(0,T,B_{2,1}^{1+\frac{d}{2}})}\\&\le c(1-e^{-(8k+2\gamma)T})\delta M\;.
\end{align*}
Collecting the estimates above and inserting them in (\ref{3.133}), we obtain 
\begin{align*}
\|w^{(n+1)}\|_{L^{1}(0,T,B_{2,1}^{\frac{d}{2}})}\le \frac{\delta}{2}+c(1-e^{-(8k+2\gamma)T})\delta +c(1-e^{-(8k+2\gamma)T})\delta M\;.
\end{align*}
Choosing T sufficiently small such that
$c(1-e^{-(8k+2\gamma)T})\le\min(\frac{1}{4M},\frac{1}{4})\;,$ we get
\begin{align*}
\|w^{(n+1)}\|_{L^{1}(0,T,B_{2,1}^{\frac{d}{2}})}\le \frac{\delta}{2}+\frac{\delta}{4}+\frac{\delta}{4}=\delta\;.
\end{align*}
These uniform bounds allow us to extract a weakly convergent subsequence. That is there is $(u,w)\in Y$ such that a subsequence  of $(u^{n},w^{n})$ (still denoted by $(u^{n},w^{n})$) satisfies 
\begin{align*}
&u^{n}\overset{*}{\rightharpoonup} u \quad\text{in}\quad L^{\infty}(0,T,B_{2,1}^{1+\frac{d}{2}-2\alpha})\;,\\&w^{n}\overset{*}{\rightharpoonup} w \quad\text{in}\quad L^{\infty}(0,T,B_{2,1}^{\frac{d}{2}})\;.
\end{align*}
In order to show that $(u,w)$ is a weak solution of (\ref{2.1a}) we need to further extract a subsequence which converges strongly to $(u,w).$ We use the Aubin-Lions Lemma. We can show by making use of the equation (\ref{2.6}) that $(\partial_tu^n,\partial_tw^n)$ is uniformly bounded in 
\begin{align*}
&\partial_tu^{n}\in L^{1}(0,T,B_{2,1}^{1+\frac{d}{2}-2\alpha})\cap L^{2}(0,T,B_{2,1}^{1+\frac{d}{2}-3\alpha})\;,\\&\partial_tw^{n}\in L^{1}(0,T,B_{2,1}^{\frac{d}{2}})\cap L^{2}(0,T,B_{2,1}^{\frac{d}{2}})\;.
\end{align*}
Since we are in this case in the whole space $\mathbb{R}^d,$ we need to combine Cantor's diagonal process with the Aubin-Lions Lemma to show that a subsequence of a weakly convergent subsequence, still denoted by $(u^n,w^n),$ has the following strongly convergent property
\begin{align*}
(u^{n},w^n)\longrightarrow(u,w) \quad \text{in}\quad L^{2}(0,T,B_{2,1}^{1+\frac{d}{2} -\gamma}(Q))\;,
\end{align*}
where $\alpha\le \gamma\le 3\alpha$ and $Q\subset\mathbb{R}^d$ is a compact subset. This strong convergence property would allow us to show that $(u,w)$ is indeed a weak solution of (\ref{2.1a}). This completes the proof  for the existence 
part of Theorem \ref{thm2}.
\end{proof}

\subsection{Uniqueness of weak solutions}
\begin{proof}[Proof]
Assume that $(u^{(1)},w^{(1)})$ and $(u^{(2)},w^{(2)})$ are two solutions of (\ref{2.1a}) in the regularity class in (\ref{2.2}) and (\ref{2.3}).
Their difference $(\widetilde{u},\widetilde{w})$ with $$\widetilde{u}=u^{(2)}-u^{(1)}\quad\text{and} \quad\widetilde{w}=w^{(2)}-w^{(1)}$$ satisfies 
\begin{align}
 \left\{
 \begin{array}{ll}
 \partial_t\widetilde{u}+(\nu+k)(-\Delta)^{\alpha}\widetilde{u}=-\mathbb{P}(u^{(2)}\cdot\nabla\widetilde{u}+\widetilde{u}\cdot\nabla u^{(1)})+2k\nabla\times\widetilde{w},\\\\
   \partial_t\widetilde{w}=-(4k+\gamma)\widetilde{w}-2k\nabla\times\widetilde{u}-u^{(2)}\cdot\nabla\widetilde{w}-\widetilde{u}\cdot \nabla w^{(1)},\\\\
\nabla\cdot\widetilde{u}=0,\\\\
\widetilde{u}(x,0)=0,\;\;\widetilde{w}(x,0)=0~.
\end{array}\label{3.4}
\right.
\end{align}
We estimate the difference $(\widetilde{u},\widetilde{w})$ in $L^2(\mathbb{R}^d)$.
Dotting (\ref{3.4}) by $(\widetilde{u},\widetilde{w})$ and applying the divergence-free condition, we find
\begin{align*}
\frac{1}{2}&\frac{d}{dt}\Big(\|\widetilde{u}\|^2_{L^2}+\|\widetilde{w}\|^2_{L^2}\Big)+(\nu+k)\|\Lambda^{\alpha}\widetilde{u}\|^2_{L^2}+(4k+\gamma)\|\widetilde{w}\|^2_{L^2}\\&=-\int u^{(2)}\cdot\nabla\widetilde{u}\cdot\widetilde{u}\,dx-\int\widetilde{u}\cdot\nabla u^{(1)}\cdot\widetilde{u}\,dx\\&\quad-\int u^{(2)}\cdot\nabla\widetilde{w}\cdot\widetilde{w}\,dx-\int\widetilde{u}\cdot\nabla w^{(1)}\cdot\widetilde{w}\,dx\\&=L_1+L_2+L_3+L_4~,
\end{align*}
where 
\begin{align*}
&L_1=-\int u^{(2)}\cdot\nabla\widetilde{u}\cdot\widetilde{u}\,dx\;,\\&L_2=-\int\widetilde{u}\cdot\nabla u^{(1)}\cdot\widetilde{u}\,dx\;,\\&L_3=-\int u^{(2)}\cdot\nabla\widetilde{w}\cdot\widetilde{w}\,dx\;,\\&L_4=-\int\widetilde{u}\cdot\nabla w^{(1)}\cdot\widetilde{w}\,dx~.
\end{align*}
Due to $\nabla\cdot u^{(2)}=0,$  we find $L_1=L_3=0$ after integration by parts. As in (\ref{ddff}), 
\begin{align*}
|L_2| \le c \|u^{(1)}\|_{B_{2,1}^{1+\frac{d}{2}}}\|\widetilde{u}\|^2_{L^2}.
\end{align*}
To bound $L_4\;,$ We set 
\begin{align*}
\frac{1}{p}=\frac{1}{2}-\frac{\alpha}{d}\;,\quad\frac{1}{q}=\frac{\alpha}{d}\;\;(or\; \frac{d}{q}=\alpha)~.
\end{align*}
By H\"older's inequality,
\begin{align*}
|L_4|&=|-\int\widetilde{u}\cdot\nabla w^{(1)}\cdot\widetilde{w}\,dx|\\&\le\|\widetilde{w}\|_{L^2}\|\nabla w^{(1)}\|_{L^q} \|\widetilde{u}\|_{L^p}\\&\le \sum_{j\ge-1}\|\Delta_j\nabla w^{(1)}\|_{L^q}\|\widetilde{w}\|_{L^2}\|\widetilde{u}\|_{L^p}\\&\le c\sum_{j\ge-1}2^j2^{dj(\frac{1}{2}-\frac{1}{q})}\|\Delta_j w^{(1)}\|_{L^2}\|\widetilde{w}\|_{L^2}\|\widetilde{u}\|_{L^p}\\&=\sum_{j\ge-1}2^j2^{\frac{dj}{2}-\frac{d}{q}j}\|\Delta_j w^{(1)}\|_{L^2}\|\widetilde{w}\|_{L^2}\|\widetilde{u}\|_{L^p}\\&\underset{\text{since}\,\alpha\ge 1}{\le} c\;\;\sum_{j\ge-1}2^{\frac{dj}{2}}\|\Delta_j w^{(1)}\|_{L^2}\|\widetilde{w}\|_{L^2}\|\widetilde{u}\|_{L^p}\\&\le c\, \|w^{(1)}\|_{B_{2,1}^{\frac{d}{2}}}\|\widetilde{w}\|_{L^2}\|\Lambda^{\alpha}\widetilde{u}\|_{L^2}\\&\le \frac{(\nu+k)}{2}\|\Lambda^{\alpha}\widetilde{u}\|_{L^2}^2+c\|w^{(1)}\|^2_{B_{2,1}^{\frac{d}{2}}}\|\widetilde{w}\|^2_{L^2}~,
\end{align*}
Where in the last inequality we make use of
\begin{align*}
\|\widetilde{u}\|_{L^p}\le c\;\|\Lambda^{\alpha}\widetilde{u}\|_{L^2}~.
\end{align*}
Combining the estimates leads to
\begin{align}
\frac{d}{dt}\Big(\|\widetilde{u}\|^2_{L^2}&+\|\widetilde{w}\|^2_{L^2}\Big)+(\nu+k)\|\Lambda^{\alpha}\widetilde{u}\|^2_{L^2}+(8k+2\gamma)\|\widetilde{w}\|^2_{L^2}\nonumber\\&\le\Big( c \|u^{(1)}\|_{B_{2,1}^{1+\frac{d}{2}}}+c\,\|w^{(1)}\|^2_{B_{2,1}^{\frac{d}{2}}}\Big)\Big(\|\widetilde{u}\|^2_{L^2}+\|\widetilde{w}\|^2_{L^2}\Big)~.\label{4.55}
\end{align}
Since $u^{(1)}\in L^1(0,T,B_{2,1}^{1+\frac{d}{2}})$ and $w^{(1)}\in L^1(0,T,B_{2,1}^{\frac{d}{2}})\cap L^{\infty}(0,T, B_{2,1}^{\frac{d}{2}})\;,$
\begin{align*}
\int_0^T\|u^{(1)}(t)\|_{B_{2,1}^{1+\frac{d}{2}}}dt<\infty \quad\quad\text{and}\quad\quad \int_0^T\|w^{(1)}(t)\|^2_{B_{2,1}^{\frac{d}{2}}}dt\le T \|w^{(1)}(t)\|^2_{L^{\infty}(0,T,B_{2,1}^{\frac{d}{2}})}\;.
\end{align*}
Applying Gronwall's inequality to (\ref{4.55}) yields
\begin{align*}
\|\widetilde{u}\|_{L^2}=\|\widetilde{w}\|_{L^2}=0\;,
\end{align*}
which leads to the desired uniqueness. This completes the proof of the uniqueness part of Theorem \ref{thm2}.
\end{proof}

\end{document}